\documentclass[dvips,oneside,noinfoline]{article}

\RequirePackage[OT1]{fontenc}
\RequirePackage[ps]{imsart20040929}
\RequirePackage[psamsfonts]{amssymb}
\usepackage{amsthm}

\RequirePackage{hyperref}

\pubyear{2004}
\volume{1}
\paperno{2}
\firstpage{20}
\lastpage{71}
\doi{10.1214/154957804100000024}

\begin{document}

\begin{frontmatter}

\title{General state space Markov chains
and MCMC algorithms\protect\thanks{This is an original survey
paper.}}
\runtitle{Markov chains
and MCMC algorithms}

\begin{aug}
\author{\fnms{Gareth O.} \snm{Roberts}\corref{}\ead[label=e1]{g.o.roberts@lancaster.ac.uk}}
\address{Department of Mathematics and Statistics,
Fylde College, Lancaster University, Lancaster, LA1 4YF, England\\
\printead{e1}}
\author{\fnms{Jeffrey S.} \snm{Rosenthal}\ead[label=e2]{jeff@math.toronto.edu}\protect\thanks{Web: \url{http://probability.ca/jeff/}.
Supported in part by NSERC of Canada.}}
\address{Department of Statistics, University of Toronto,
Toronto, Ontario, Canada  M5S 3G3\\
\printead{e2}}
\end{aug}


\runauthor{G.O. Roberts, J.S. Rosenthal}

\begin{abstract}
This paper surveys various results about Markov chains on general
(non-countable) state spaces.  It begins with an introduction to Markov
chain Monte Carlo (MCMC) algorithms, which provide the motivation and
context for the theory which follows.  Then, sufficient conditions for
geometric and uniform ergodicity are presented, along with quantitative
bounds on the rate of convergence to stationarity.  Many of these results
are proved using direct coupling constructions based on minorisation
and drift conditions.  Necessary and sufficient conditions for Central
Limit Theorems (CLTs) are also presented, in some cases proved via
the Poisson Equation or direct regeneration constructions.  Finally,
optimal scaling and weak convergence results for Metropolis-Hastings
algorithms are discussed.  None of the results presented is new, though
many of the proofs are.  We also describe some Open Problems.
\end{abstract}
\received{\smonth{3} \syear{2004}}



\end{frontmatter}
\def\suml{\sum\limits}
\def\prodl{\prod\limits}
\def\intl{\int\limits}
\def\supl{\sup\limits}
\def\infl{\inf\limits}
\def\liml{\lim\limits}
\def\maxl{\max\limits}
\def\minl{\min\limits}
\def\real{{\Re e \, }}
\def\tr{{\rm tr \, }}
\def\vecx{\vec{x}}
\def\vecy{\vec{y}}
\def\zeroone{\{0,1\}}
\def\overn{{1 \over n}}
\def\darkbox{{\vrule height8pt width6pt depth0pt}}
\def\half { { 1 \over 2 } }
\def\un{\underbar}
\def\const {(\hbox{const})}
\def\blank {\underbar{ \ \ }}
\def\vark {{ \| \mu_k - \pi \| }}

\def\bX{{\bf X}}
\def\bY{{\bf Y}}
\def\bx{{\bf x}}
\def\by{{\bf y}}
\def\bz{{\bf z}}
\def\bZ{{\bf Z}}
\def\bw{{\bf w}}
\def\bW{{\bf W}}
\def\ba{{\bf a}}
\def\btheta{{\bf \theta}}
\def\IR{{\bf R}}
\def\IC{{\bf C}}
\def\IZ{{\bf Z}}
\def\IN{{\bf N}}
\def\X{\mathcal{X}}
\def\L{{\cal L}}
\def\Y{{\cal Y}}
\def\U{{\cal U}}
\def\F{{\cal F}}
\def\A{{\cal A}}
\def\Zn{{\IZ/(n)}}
\def\E{{\bf E}}
\def\P{{\bf P}}
\def\Var{{\bf Var}}
\def\Cov{{\bf Cov}}
\def\Exp{{\bf Exp}}
\def\one{{\bf 1}}
\def\point{\medskip{$\bullet$}}
\def\grad{\nabla}

\def\Pb{\hat{P}}
\def\Qb{\hat{Q}}

\def\Pbar{\overline{P}}
\def\hbar{\overline{h}}
\def\pibar{\overline{\pi}}

\def\state{\mathcal{X}}
\def\be{\begin{equation}}
\def\ee{\end{equation}}
\def\beq{\begin{eqnarray*}}
\def\eeq{\end{eqnarray*}}
\def\sc{\setcounter}
\def\Corr{\mathop{\mathrm{Corr}}\nolimits}
\def\bye{\end{document}}

\newtheorem{thm}{Theorem}
\newtheorem{Proposition}{Proposition}
\newtheorem{Lemma}{Lemma}
\newtheorem{Corollary}{Corollary}
\newtheorem{opennew}{Open Problem \#}
\theoremstyle{definition}
\newtheorem*{defn}{Definition}
\newtheorem{Claim}{Claim}
\newtheorem{Fact}{Fact}
\newtheorem{ex}{Example}
\newtheorem*{remark}{Remark}
\newtheorem*{Running Example}{Running Example}
\newtheorem*{Running Example, Continued}{Running Example, Continued}
\newtheorem*{PfC1}{Proof of Claim~\ref{Cl1}}
\newtheorem*{PfC2}{Proof of Claim~\ref{Cl2}}
\newtheorem*{PfT25}{Proof of Theorem~\ref{reversibleclt}}
\newtheorem*{PfT27}{Proof of Theorem~\ref{VCLTthm}}
\newtheorem*{proof1}{Proof}
\newtheorem*{PfL17}{Proof of Lemma~\ref{petitelemma2}}

\section{Introduction}
\label{S1}

Markov chain Monte Carlo (MCMC) algorithms -- such as the
Metropolis-Hastings algorithm (\cite{metropolis}, \cite{hastings})
and the Gibbs sampler (e.g.\ Geman and Geman~\cite{geman}; Gelfand and
Smith~\cite{gelfand}) -- have become extremely popular in statistics, as a
way of approximately sampling from complicated probability distributions
in high dimensions (see for example the reviews \cite{tierney},
\cite{smithrob}, \cite{gilksbook}, \cite{bruns}).  Most dramatically,
the existence of MCMC algorithms has transformed Bayesian inference,
by allowing practitioners to sample from posterior distributions of
complicated statistical models.

In addition to their importance to applications in statistics and
other subjects, these algorithms also raise numerous questions
related to probability theory and the mathematics of Markov chains.
In particular, MCMC algorithms involve Markov chains $\{X_n\}$ having
a (complicated) stationary distribution $\pi(\cdot)$, for which it is
important to understand as precisely as possible the nature and speed
of the convergence of the law of $X_n$ to $\pi(\cdot)$ as $n$ increases.

This paper attempts to explain and summarise MCMC algorithms and the
probability theory questions that they generate.  After introducing
the algorithms (Section~\ref{S2}), we discuss various important theoretical
questions related to them.  In
Section~\ref{S3} we present various convergence
rate results commonly used in MCMC.  Most of these are proved in
Section~\ref{S4}, using direct coupling arguments and thereby avoiding many of
the analytic technicalities of previous proofs.  We consider MCMC central
limit theorems in Section~\ref{S5}, and optimal scaling and weak convergence
results in Section~\ref{S6}.  Numerous references to the MCMC literature are
given throughout.
We also describe some Open Problems.


\subsection{The problem}
\label{S1.1}

The problem addressed by MCMC algorithms is the following.
We're given a density function $\pi_u$, on some
state space $\X$, which is possibly unnormalised but at least
satisfies $0 < \int_\X \pi_u < \infty$.
(Typically $\X$ is an open subset of $\IR^d$, and
the densities are taken with respect to Lebesgue measure, though other
settings -- including discrete state spaces -- are also possible.)  This
density gives rise to a probability measure $\pi(\cdot)$ on $\X$, by
\be
\pi(A)  =  {\int_{A} \pi_u(x) dx \over \int_{\state} \pi_u (x) dx} .
\label{density}
\label{e1}
\ee

We want to (say)
estimate expectations of functions $f:\X\to\IR$
with respect to $\pi(\cdot)$, i.e.\ we want to estimate
\be
\pi(f)  =  \E_{\pi}[f(X)]
 =  {\int_{\state} f(x) \pi_u (x) dx \over \int_{\state} \pi_u (x) dx} .
\label{pifdef}
\label{e2}
\ee
If $\state$ is high-dimensional, and $\pi_u$ is a
complicated function, then direct integration (either analytic or numerical)
of the integrals in~(\ref{pifdef}) is infeasible.

The classical Monte Carlo solution to this problem is to simulate
i.i.d.\ random variables $Z_1,Z_2,\ldots,Z_N \sim \pi(\cdot)$, and
then estimate $\pi(f)$ by
\be
\hat\pi(f)  =  (1/N) \sum_{i=1}^N f(Z_i)  .
\label{classicalmc}
\label{e3}
\ee
This gives an unbiased estimate, having standard deviation
of order $O(1/\sqrt{N})$.  Furthermore, if $\pi(f^2)<\infty$, then
by the classical Central Limit Theorem, the error $\hat\pi(f)
- \pi(f)$ will have a limiting normal distribution, which is also useful.
The problem, however, is that if $\pi_u$ is complicated, then it is very
difficult to directly simulate i.i.d.\ random variables from $\pi(\cdot)$.

The Markov chain Monte Carlo (MCMC) solution is to instead
construct a {\it Markov chain} on $\state$ which is easily run on a
computer, and which has $\pi(\cdot)$ as a stationary distribution.
That is, we want to define easily-simulated Markov chain
transition probabilities $P(x,dy)$ for $x,y\in\state$, such that
\be
\int_{x \in \state} \pi(dx) \, P(x,dy)  =  \pi(dy) .
\label{stationaryeqn}
\label{e4}
\ee
Then hopefully (see Subsection~\ref{S3.2}),
if we run the Markov chain for a long time (started
from anywhere), then for large $n$ the distribution of $X_n$ will be
approximately stationary: $\L(X_n) \approx \pi(\cdot)$.  We can then
(say) set $Z_1 = X_n$, and then
restart and rerun the Markov chain to obtain $Z_2$, $Z_3$, etc., and
then do estimates as in~(\ref{classicalmc}).

It may seem at first to be even more difficult to find such a Markov
chain, then to estimate $\pi(f)$ directly.  However, we shall see in
the next section that constructing (and running) such Markov chains is
often surprisingly straightforward.

\begin{remark}
In the practical use of MCMC, rather than start a fresh Markov chain for
each new sample, often an entire tail of the Markov chain run $\{X_n\}$
is used to create an estimate such as $(N-B)^{-1} \sum_{i=B+1}^N f(X_i)$,
where the {\it burn-in} value $B$ is hopefully
chosen large enough that $\L(X_B) \approx \pi(\cdot)$.
In that case the different $f(X_i)$ are not independent, but the estimate
can be computed more efficiently.  Since many of the mathematical issues
which arise are similar in either implementation, we largely ignore
this modification herein.
\end{remark}

\begin{remark}
MCMC is, of course, not the only way to sample or estimate from complicated
probability distributions.  Other possible sampling algorithms include
``rejection sampling'' and ``importance sampling'', not reviewed here;
but these alternative algorithms only work well in certain particular
cases and are not as widely applicable as MCMC algorithms.
\end{remark}

\subsection{Motivation: Bayesian Statistics Computations}
\label{S1.2}

While MCMC algorithms are used in many fields (statistical physics,
computer science), their most widespread application is in Bayesian
statistical inference.

Let $L({\bf y}|{\bf \theta })$ be the likelihood function (i.e., density
of data $\by$ given unknown parameters $\theta$) of a statistical model,
for ${\bf \theta}\in \X$.  (Usually $\X \subseteq {\bf R}^d$.)
Let the ``prior'' density of ${\bf \theta}$ be $p({\bf \theta})$.
Then the ``posterior'' distribution of ${\bf \theta}$ given ${\bf y}$
is the density which is proportional to
\[
\pi_u({\bf \theta})
 \equiv
L({\bf y} \mid {\bf \theta }) \, p({\bf \theta}) .
\]
(Of course, the normalisation constant is simply the density for the
data $\by$, though that constant may be impossible to compute.)
The ``posterior mean'' of any functional $f$ is then given by:
\[
\pi(f)  =
{\int_{\state} f(x) \pi_u (x) dx \over \int_{\state} \pi_u (x) dx}  .
\]

For this reason, Bayesians are anxious (even desperate!) to estimate
such $\pi(f)$.  Good estimates allow Bayesian inference can be used
to estimate a wide variety of parameters, probabilities, means, etc.
MCMC has proven to be extremely helpful for such Bayesian estimates, and
MCMC is now extremely widely used in the Bayesian statistical community.

\section{Constructing MCMC Algorithms}
\label{S2}

We see from the above that an MCMC algorithm requires, given a
probability distribution $\pi(\cdot)$ on a state space $\X$, a Markov
chain on $\state$ which is easily run on a computer, and which has
$\pi(\cdot)$ as its stationary distribution as
in~(\ref{stationaryeqn}).
This section explains how such Markov chains are constructed.  It thus
provides motivation and context for the theory which follows; however,
for the reader interested purely in the mathematical results, this
section can be omitted with little loss of continuity.

A key notion is {\it reversibility}, as follows.

\begin{defn}
A Markov chain on a state space $\X$ is {\it reversible} with respect to a
probability distribution $\pi(\cdot)$ on $\X$, if
\[
\pi(dx) \, P(x,dy)  =  \pi(dy) \, P(y,dx), \qquad x,y\in\state  .
\]

A very important property of reversibility is the following.
\end{defn}

\begin{Proposition}
\label{revprop}\label{P1}
If Markov chain is reversible with respect to $\pi(\cdot)$, then
$\pi(\cdot)$ is stationary for the chain.
\end{Proposition}

\begin{proof}
We compute that
\[
\int_{x \in \state} \pi(dx) \, P(x,dy)
 =  \int_{x \in \state} \pi(dy) \, P(y,dx)
 =  \pi(dy) \, \int_{x \in \state} P(y,dx)
=  \pi(dy)  .
\]
\end{proof}

We see from this lemma that, when constructing an MCMC algorithm, it
suffices to create a Markov chain which is easily run, and which is
{\it reversible} with respect to $\pi(\cdot)$.  The simplest way to do
so is to use the Metropolis-Hastings algorithm, as we now discuss.

\subsection{The Metropolis-Hastings Algorithm}
\label{S2.1}

Suppose again that $\pi(\cdot)$ has a (possibly unnormalised) density
$\pi_u$, as in~(\ref{density}).  Let $Q(x,\cdot)$ be essentially any
other Markov chain, whose transitions also have a (possibly unnormalised)
density, i.e.\ $Q(x,dy) \propto q(x,y) \, dy$.

The Metropolis-Hastings algorithm proceeds as follows. First choose
some~$X_0$.  Then, given $X_n$,
generate a {\it proposal}\/ $Y_{n+1}$ from $Q(X_n, \cdot)$.
Also flip an independent coin, whose probability of heads equals
$\alpha(X_n, Y_{n+1})$, where
\[
\alpha(x, y)=\min\bigg[1, \ {\pi_u(y) \, q(y, x) \over \pi_u(x) \, q(x, y)}
\bigg] .
\]
(To avoid ambiguity, we set $\alpha(x,y)=1$ whenever $\pi(x) \, q(x, y) = 0$.)
Then, if the coin is heads, ``accept'' the proposal by
setting $X_{n+1} = Y_{n+1}$; if the coin is tails then
``reject'' the proposal by setting $X_{n+1} = X_n$.
Replace $n$ by $n+1$ and repeat.

The reason for the unusual formula for $\alpha(x,y)$ is the following:

\begin{Proposition}
\label{P2}
The Metropolis-Hastings algorithm (as described above) produces a Markov
chain $\{X_n\}$ which is reversible with respect to $\pi(\cdot)$.
\end{Proposition}
\begin{proof}
We need to show
\[
\pi(dx) \, P(x,dy) = \pi(dy) \, P(y,dx) .
\]
It suffices to assume $x \not= y$ (since if $x=y$ then the equation is
trivial).  But for $x\not=y$, setting $c = \int_\X \pi_u(x) \, dx$,
\begin{eqnarray*}
&&\pi(dx) \, P(x,dy)
= [c^{-1} \pi_u(x) \, dx] \, [q(x,y) \, \alpha(x,y) \, dy]\\
&&\qquad = c^{-1} \pi_u(x) \, q(x,y)
\, \min\bigg[1, \ {\pi_u(y) q(y, x) \over \pi_u(x) q(x, y)} \bigg]
\, dx \, dy\\
&&\qquad = c^{-1} \min[ \pi_u(x) \, q(x,y), \ \pi_u(y) q(y, x) ]\, dx \, dy ,
\end{eqnarray*}
which is symmetric in $x$ and $y$.
\end{proof}

To run the Metropolis-Hastings algorithm on a computer, we just need
to be able to run the proposal chain $Q(x,\cdot)$ (which is easy, for
appropriate choices of $Q$), and then do the accept/reject step (which is
easy, provided we can easily compute the densities at individual points).
Thus, running the algorithm is quite feasible.
Furthermore we need to compute only {\it ratios} of densities [e.g.\
$\pi_u(y) \, / \, \pi_u(x)$], so we don't require the {\it normalising
constants} $c=\int_\X \pi_u(x) dx$.

However, this algorithm in turn suggests further questions.  Most
obviously, how should we choose the proposal distributions $Q(x,\cdot)$?
In addition, once $Q(x,\cdot)$ is chosen, then will we really have
$\L(X_n) \approx \pi(\cdot)$ for large enough $n$?  How large is large
enough?  We will return to these questions below.

Regarding the first question, there are many different classes of ways
of choosing the proposal density, such as:

\point
{\bf Symmetric Metropolis Algorithm.}
Here
$
q(x, y) = q(y,x)
$, and
the acceptance probability simplifies to
\[
\alpha (x, y) = \min\bigg[1, \ {\pi_u(y)\over \pi_u(x)}\bigg]
\]

\point
{\bf Random walk Metropolis-Hastings.}
Here $q(x, y) = q(y-x)$.  For example, perhaps
$Q(x,\cdot) = N(x,\sigma^2)$, or
$Q(x,\cdot) = {\rm Uniform}(x-1, \, x+1)$.

\point
{\bf Independence sampler.}
Here $q(x, y) = q(y)$, i.e.\ $Q(x,\cdot)$ does not depend on~$x$.

\point
{\bf Langevin algorithm.}
Here the proposal is generated by
\[
Y_{n+1} \sim N(X_n + (\delta/2) \, \nabla \log \pi (X_n), \ \delta) ,
\]
for some (small) $\delta>0$.
(This is motivated by a discrete approximation to a
Langevin diffusion processes.)

\medskip

More about optimal choices of proposal distributions will be discussed in
a later section, as will the second question about time to stationarity
(i.e.\ how large does $n$ need to be).

\subsection{Combining Chains}
\label{S2.2}

If $P_1$ and $P_2$ are two different chains, each having stationary
distribution $\pi(\cdot)$, then the new chain $P_1 P_2$ also has
stationary distribution $\pi(\cdot)$.

Thus, it is perfectly acceptable, and quite common (see
e.g.\ Tierney~\cite{tierney}
and \cite{hybrid}), to make new MCMC algorithms
out of old ones, by specifying that the new algorithm applies first
the chain $P_1$, then the chain $P_2$, then the chain $P_1$ again, etc.
(And, more generally, it is possible to combine many different chains
in this manner.)

Note that, even if each of $P_1$ and $P_2$ are reversible, the combined
chain $P_1 P_2$ will in general {\it not} be reversible.  It is for
this reason that it is important, when studying MCMC, to allow for
non-reversible chains as well.

\subsection{The Gibbs Sampler}
\label{S2.3}

The Gibbs sampler is also known as the ``heat bath'' algorithm,
or as ``Glauber dynamics''.
Suppose again that $\pi_u(\cdot)$ is $d$-dimensional density, with
${\state}$ an open subset of ${\bf R}^d$, and
write ${\bf x}= (x_1,\ldots ,x_d)$.

The {\bf $i^{\rm th}$ component Gibbs sampler} is defined such that $P_i$
leaves all components besides $i$ unchanged, and replaces the $i^{\rm th}$
component by a draw from the full conditional distribution of $\pi(\cdot)$
conditional on \un{all} the other components.

More formally, let
\[
S_{x,i,a,b}  =  \{y \in \X; \ y_j=x_j \ {\rm for} \ j\not=i, \
{\rm and} \ a \le y_i \le b \} .
\]
Then
\[
P_i(x, \, S_{x,i,a,b})  =
{ \int_a^b \pi_u(x_1,\ldots,x_{i-1},t,x_{i+1},\ldots,x_n) \, dt
\over
\int_{-\infty}^\infty \pi_u(x_1,\ldots,x_{i-1},t,x_{i+1},\ldots,x_n) \, dt
}  ,
\quad a \le b  .
\]

It follows immediately (from direct computation, or from the definition
of conditional density), that $P_i$, is reversible with respect to
$\pi(\cdot)$.  (In fact, $P_i$ may be regarded as a special case of a
Metropolis-Hastings algorithm, with $\alpha(x,y) \equiv 1$.)
Hence, $P_i$ has $\pi(\cdot)$ as a stationary distribution.

We then construct the full Gibbs sampler out of the various $P_i$, by
combining them (as in the previous subsection) in one of two ways:

\medskip\point
{\bf The deterministic-scan Gibbs sampler} is
\[
P  =  P_1 P_2 \ldots P_d  .
\]
That is, it performs the $d$ different Gibbs sampler components, in
sequential order.

\medskip\point
{\bf The random-scan Gibbs sampler} is
$$
P  =  {1 \over d} \sum_{i=1}^d P_i .
$$
That is, it does one of the $d$ different Gibbs sampler components, chosen
uniformly at random.

\medskip
Either version produces an MCMC algorithm having $\pi(\cdot)$ as its
stationary distribution.  The output of a Gibbs sampler is thus
a ``zig-zag pattern'', where the components get updated one at a time.
(Also, the random-scan Gibbs sampler is reversible, while the
deterministic-scan Gibbs sampler usually is not.)

\subsection{Detailed Bayesian Example: Variance Components Model}
\label{S2.4}

We close this section by presenting a typical example of a target
density $\pi_u$ that arises in Bayesian statistics, in an effort to
illustrate the problems and issues which arise.

The model involves fixed constant $\mu_0$ and positive constants
$a_1,b_1,a_2,b_2$,
and $\sigma_0^2$.  It involves three hyperparameters, $\sigma_\theta^2$,
$\sigma_e^2$, and $\mu$, each having priors based
upon these constants as follows:
\ \ $\sigma_\theta^2 \sim IG(a_1,b_1)$;
\ \ $\sigma_e^2 \sim IG(a_2,b_2)$;
\ \ and $\mu \sim N(\mu_0, \, \sigma_0^2)$.
It involves $K$ further parameters $\theta_1,\theta_2,\ldots,\theta_K$,
conditionally independent given the above hyperparameters,
with $\theta_i \sim N(\mu,\sigma_\theta^2)$.
In terms of these parameters, the data $\{Y_{ij}\}$ $(1 \le i \le K, \ 1
\le j \le J)$ are assumed to be distributed as $Y_{ij} \sim N(\theta_i,
\sigma_e^2)$, conditionally independently given the parameters.
A graphical representation of the model is as follows:

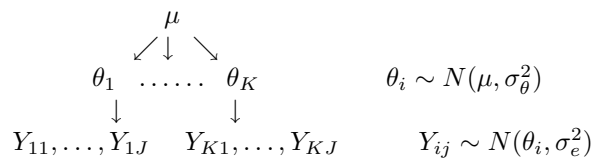
\def\mspc{ \ \ }
\[
\begin{array}{c@{\qquad\quad}c}
\hskip -0.1cm
\mu  & \cr
\swarrow \mspc \hskip -0.2cm
 \downarrow \mspc \searrow & \cr
\theta_1 \mspc \ldots \ldots \mspc \theta_K &
\hskip -1.1cm\theta_i \sim N(\mu, \sigma_\theta^2) \cr
\downarrow \mspc \mspc \mspc
 \mspc \mspc \mspc
\downarrow \cr Y_{11},\ldots,Y_{1J} \ \  \ \
Y_{K1},\ldots,Y_{KJ} & Y_{ij} \sim N(\theta_i,
\sigma_e^2) 
\end{array}
\]

The Bayesian paradigm then involves conditioning on the values of the
data $\{Y_{ij}\}$, and considering the joint distribution of all $K+3$
parameters given this data.  That is, we are interested in the distribution
$$
\pi(\cdot) = {\cal L}(\sigma_\theta^2, \sigma_e^2, \mu,
\theta_1,\ldots,\theta_K \mid \{Y_{ij}\})  ,
$$
defined on the state space $\X = (0,\infty)^2 \times {\bf R}^{K+1}$.
We would like to sample from this distribution $\pi(\cdot)$.
We compute that this distribution's unnormalised density is given by
\beq
&&\pi_u(\sigma^2_\theta, \sigma^2_e, \mu, \theta_1,\ldots,\theta_K)
\propto\\[6pt]
&& e^{-b_1/\sigma^2_\theta} {\sigma^2_\theta}^{-a_1-1}
e^{-b_2/\sigma^2_e} {\sigma^2_e}^{-a_2-1}
e^{-(\mu-\mu_0)^2 / 2\sigma_0^2}\\[6pt]
&&{}\times \prod_{i=1}^K [e^{-(\theta_i-\mu)^2 / 2 \sigma^2_\theta} /
\sigma_\theta]
\times \prod_{i=1}^K \prod_{j=1}^J [e^{-(Y_{ij}-\theta_i)^2 / 2 \sigma^2_e}
/ \sigma_e].
\eeq\vspace*{-28pt}\eject

\noindent This
is a very typical target density for MCMC in statistics, in that it
is high-dimensional ($K+3$), its formula is messy and irregular, it is
positive throughout $\X$, and it is
larger in ``center'' of $\X$ and smaller in ``tails'' of~$\X$.

We now consider constructing MCMC algorithms
to sample from the target density $\pi_u$.  We begin with the
Gibbs sampler.  To run a Gibbs sampler, we require the full conditionals
distributions, computed (without difficulty since they are all
one-dimensional) to be as follows:

\def\st2{{\sigma_\theta^2}}
\def\se2{{\sigma_e^2}}
\def\s02{\sigma_0^2}
\def\thetas{\theta_1,\ldots,\theta_K}
\def\Ybar{\overline{Y}}
\def\overJ{{1 \over J}}

\beq
&&\hspace*{-7mm}\L(\st2 \mid \mu,\se2,\thetas,Y_{ij})  =
IG \bigg( a_1 + \half K, \ b_1
+ \half \sum_i (\theta_i - \mu)^2 \bigg) ;\\
&&\hspace*{-7mm}\L(\se2 \mid \mu,\st2,\thetas,Y_{ij})  =
IG \bigg( a_2 + \half KJ, \ b_2 + \half \sum_{i,j} (Y_{ij} - \theta_i)^2
           \bigg) ;\\
&&\hspace*{-7mm}\L(\mu \mid \st2,\se2,\thetas,Y_{ij})  =
N \bigg( { \st2 \mu_0 + \sigma_0^2 \sum_i \theta_i \over \st2
+ K\sigma_0^2 }, \
{ \st2 \sigma_0^2 \over \st2 + K\sigma_0^2 }\bigg) ;\\
&&\hspace*{-7mm}\L(\theta_i \mid
\mu,\st2,\se2,\theta_1,\ldots,\theta_{i-1},\theta_{i+1},\ldots,
\theta_K,Y_{ij})  =
N \bigg( { J \st2 \overline{Y}_i + \se2 \mu \over J \st2 + \se2 }, \
{ \st2 \se2 \over J \st2 + \se2 }\bigg),
\eeq
where $\Ybar_i = \overJ \sum_{j=1}^J Y_{ij}$,
and the last equation holds for $1 \le i \le K$.
The Gibbs sampler then proceeds by updating the $K+3$ variables, in turn
(either deterministic or random scan), according to the above conditional
distributions.  This is feasible since the conditional distributions
are all easily simulated (IG and N).
In fact, it appears to work well, both in practice and according to various
theoretical results; this model was one of the early statistical
applications of the Gibbs sampler by Gelfand and Smith~\cite{gelfand},
and versions of it
have been used and studied often (see e.g.\ \cite{varcomp}, \cite{mykland},
\cite{jamesstein}, \cite{cowles1}, \cite{hobert1}, \cite{hobert2}).

\medskip

Alternatively, we can run a \un{Metropolis-Hastings algorithm} for this model.
For example, we might choose a symmetric random-walk Metropolis
algorithm with proposals of the form $N(X_n, \sigma^2 I_{K+3})$ for some
$\sigma^2>0$ (say).  Then,
given $X_n$, the algorithm would proceed as follows:
\begin{enumerate}
\item[1.]
Choose $Y_{n+1} \sim N(X_n, \sigma^2 I_{K+3})$;
\item[2.]
Choose $U_{n+1} \sim {\rm Uniform}[0,1]$;
\item[3.]
If $U_{n+1} < \pi_u(Y_{n+1}) \, / \, \pi_u(X_n)$, then set $X_{n+1} =
Y_{n+1}$ (accept).  Otherwise set $X_{n+1} = X_n$ (reject).
\end{enumerate}

This MCMC algorithm also appears to work well for this model, at least
if the value of $\sigma^2$ is chosen appropriately (as discussed in
Section~\ref{S6}).  We conclude that, for such ``typical'' target
distributions~$\pi(\cdot)$, both the Gibbs sampler and appropriate
Metropolis-Hastings algorithms perform well in practice, and allow us
to sample from~$\pi(\cdot)$.

\section{Bounds on Markov Chain Convergence Times}
\label{S3}

Once we know how to construct (and run) lots of different MCMC
algorithms, other questions arise.  Most obviously, do they converge to
the distribution $\pi(\cdot)$?  And, how quickly does this convergence
take place?

To proceed,
write $P^n(x,A)$ for the $n$-step transition law of the Markov chain:
$$
P^n(x,A) = \P[X_n \in A \mid  X_0 = x]  .
$$
The main MCMC convergence questions are, is $P^n(x,A)$ ``close'' to $\pi(A)$
for large enough $n$?  And, how large is large enough?

\subsection{Total Variation Distance}
\label{S3.1}

We shall measure the distance to stationary
in terms of total variation distance, defined as follows:

\begin{defn}
The {\it total variation distance} between two probability measures
$\nu_1(\cdot)$ and $\nu_2(\cdot)$ is:
$$
\| \nu_1(\cdot) -\nu_2(\cdot) \| = \sup _A |\nu _1(A) - \nu _2 (A) |  .
$$
\end{defn}

We can then ask, is $\lim_{n\to\infty} \| P^n(x,\cdot) - \pi(\cdot) \| =
0$?  And, given $\epsilon>0$, how large must $n$ be so that
$\| P^n(x,\cdot) - \pi(\cdot) \| < \epsilon$?  We consider such
questions herein.

We first pause to note some simple properties of total variation distance.

\begin{Proposition}
\label{tvprop}\label{P3}
(a)
$\| \nu_1(\cdot) -\nu_2(\cdot) \| = \sup_{f:\X\to[0,1]} | \int f d\nu_1
- \int f d\nu_2 |$.
\hfil\break (b)
$\| \nu_1(\cdot) -\nu_2(\cdot) \| = {1 \over b-a} \sup_{f:\X\to[a,b]}
| \int f d\nu_1 - \int f d\nu_2 |$ for any $a<b$, and in particular
$\| \nu_1(\cdot) -\nu_2(\cdot) \| = \half \sup_{f:\X\to[-1,1]} | \int f d\nu_1
- \int f d\nu_2 |$.
\hfil\break (c) If $\pi(\cdot)$ is stationary for a Markov chain
kernel~$P$, then
$\| P^n(x,\cdot) - \pi(\cdot) \|$ is non-increasing in $n$, i.e.\
$\| P^{n}(x,\cdot) - \pi(\cdot) \| \le \| P^{n-1}(x,\cdot) - \pi(\cdot) \|$
for $n\in\IN$.
\hfil\break (d) More generally, letting $(\nu_i P)(A) = \int \nu_i(dx)
\, P(x,A)$, we always have
$\| (\nu_1 P)(\cdot) - (\nu_2 P)(\cdot) \|
\le \| \nu_1(\cdot) - \nu_2(\cdot) \|$.
\hfil\break (e) Let $t(n) = 2 \, \sup_{x\in\X} \|P^n(x,\cdot)-\pi(\cdot)\|$,
where $\pi(\cdot)$ is stationary.  Then $t$ is sub-multiplicative, i.e.\
$t(m+n) \le t(m) \, t(n)$ for $m,n\in\IN$.
\hfil\break (f) If $\mu(\cdot)$ and $\nu(\cdot)$ have densities $g$ and
$h$, respectively, with respect to some $\sigma$-finite measure
$\rho(\cdot)$, and $M=\max(g,h)$ and $m=\min(g,h)$, then
$$
\|\mu(\cdot) - \nu(\cdot)\|  =  \half \int_\X (M-m) \, d\rho
 =  1 - \int_\X m \, d\rho.
$$
\hfil\break (g) Given probability measures $\mu(\cdot)$ and
$\nu(\cdot)$, there are jointly defined random variables $X$ and $Y$
such that $X \sim \mu(\cdot)$, $Y \sim \nu(\cdot)$, and
$\P[X = Y] = 1 - \|\mu(\cdot) - \nu(\cdot)\|$.
\end{Proposition}\vspace*{-12pt}\eject

\begin{proof}
For (a), let $\rho(\cdot)$ be any $\sigma$-finite measure
such that $\nu_1 \ll \rho$ and $\nu_2 \ll \rho$ (e.g.\ $\rho = \nu_1 +
\nu_2$), and set $g = d\nu_1/d\rho$ and $h = d\nu_2/d\rho$.  Then
$| \int f d\nu_1 - \int f d\nu_2 | = | \int f(g-h) \, d\rho |$.  This is
maximised (over all $0 \le f \le 1$) when $f=1$ on $\{g>h\}$ and
$f=0$ on $\{h>g\}$ (or vice-versa), in which case it equals
$|\nu_1(A)-\nu_2(A)|$
for $A = \{g>h\}$ (or $\{g<h\}$), thus proving the equivalence.

Part~(b) follows very similarly to~(a), except now $f=b$ on $\{g>h\}$
and $f=a$ on $\{g<h\}$ (or vice-versa), leading to
$| \int f d\nu_1 - \int f d\nu_2 | = (b-a) \, |\nu_1(A)-\nu_2(A)|$.

For part~(c), we compute that
\beq
&&| P^{n+1}(x,A) - \pi(A) |
 =
\left| \int_{y\in \X} P^n(x,dy) P(y,A)
- \int_{y\in \X} \pi(dy) P(y,A) \right|\\
&&\qquad  =  \left| \int_{y\in \X} P^n(x,dy) f(y)
- \int_{y\in \X} \pi(dy) f(y) \right|
 \le  \|P^n(x,\cdot) - \pi(\cdot)\| ,
\eeq
where $f(y) = P(y,A)$, and where the inequality comes from part~(a).

Part~(d) follows very similarly to part~(c).

Part~(e) follows since $t(n)$ is an $L^\infty$ operator norm
of~$P^n$ (cf.\ Meyn and Tweedie~\cite{MT}, Lemma~16.1.1).
More specifically, let $\Pb(x,\cdot) = P^n(x,\cdot) - \pi(\cdot)$
and $\Qb(x,\cdot) = P^m(x,\cdot) - \pi(\cdot)$, so that
\beq
(\Pb \, \Qb f)(x)
& \equiv &
\int_{y\in\X} f(y) \int_{z\in\X} [P^n(x,dz)-\pi(dz)] \,
[\P^m(z,dy)-\pi(dy)]\\
& = & \int_{y\in\X} f(y) \, [P^{n+m}(x,dy) - \pi(dy) - \pi(dy) +
\pi(dy)]\\
& = & \int_{y\in\X} f(y) \, [P^{n+m}(x,dy) - \pi(dy)].
\eeq
Then let $f:\X\to[0,1]$, let
$g(x) = (\Qb f)(x) \equiv \int_{y\in\X} \Qb(x,dy) f(y)$, and let
$g^* = \sup_{x\in\X} |g(x)|$.
Then $g^* \le \half \, t(m)$ by part~(a).
Now, if $g^*=0$, then clearly $\Pb\Qb f = 0$.
Otherwise, we compute that
\be
2 \, \sup_{x\in\X} | (\Pb \, \Qb f)(x) |
= 2 \, g^* \, \sup_{x\in\X} | (\Pb [g/g^*])(x) |
\le t(m) \, \sup_{x\in\X} (\Pb [g/g^*])(x) | .
\label{PBineq}
\label{e5}
\ee
Since $-1 \le g/g^* \le 1$, we have
$(\Pb [g/g^*])(x) \le 2 \, \|P^n(x,\cdot)-\pi(\cdot)\|$ by part~(b), so that
$\sup_{x\in\X} (\Pb [g/g^*])(x) \le t(n)$.
The result then follows from part~(a) together with~(\ref{PBineq}).

The first equality of part~(f) follows since, as in the proof of
part~(b) with $a=-1$ and $b=1$, we have
$$
\|\mu(\cdot) - \nu(\cdot)\|
 =  \half \Big( \int_{g>h} (g-h) \, d\rho + \int_{h>g} (h-g) \, d\rho \Big)
 =  \half \int_{\X} (M-m) \, d\rho.
$$
The second equality of part~(f) then follows since $M+m=g+h$, so that
$\int_\X (M+m) \, d\rho = 2$, and hence\eject
\beq
\half \int_{\X} (M-m) \, d\rho
& = & 1 - \half \bigg( 2 - \int_\X (M-m) \, d\rho \bigg)\\
& = & 1 - \half \int_\X \big( (M+m)-(M-m) \big) \, d\rho
 =  1 - \int_\X m \, d\rho .
\eeq

For part~(g), we let $a = \int_\X m \, d\rho$,
$b = \int_\X (g-m) \, d\rho$, and
$c = \int_\X (h-m) \, d\rho$.  The statement is trivial if any of
$a,b,c$ equal zero, so assume they are all positive.
We then jointly construct random
variables $Z,U,V,I$ such that $Z$ has density $m/a$,
$U$ has density $(g-m)/b$,
$V$ has density $(h-m)/b$,
and $I$ is independent of $Z,U,V$ with $\P[I=1]=a$ and
$\P[I=0]=1-a$.  We then let $X=Y=Z$ if $I=1$, and
$X=U$ and $Y=V$ if $I=0$.  Then it is easily checked that
$X \sim \mu(\cdot)$ and $Y \sim \nu(\cdot)$.  Furthermore
$U$ and $V$ have disjoint support, so $\P[U=V] = 0$.  Then
using part~(f),
$$
\P[X=Y]  =  \P[I=1]  =  a  =  1 \, - \, \|\mu(\cdot) - \nu(\cdot)\|
 ,
$$
as claimed.
\end{proof}

\begin{remark}
Proposition~\ref{P3}(e) is false without the factor of~2.  For
example, suppose ${\cal X}=\{1, 2\}$, with $P(1, \{1\})=0.3$, $P(1,\{2\})=0.7$,
$P(2, \{1\})=0.4$, $P(2, \{2\})=0.6$, $\pi(1)={4\over {11}}$, and
$\pi(2)={7\over {11}}$.  Then $\pi(\cdot)$ is stationary, and
$\sup_{x\in {\cal X}} \|P(x, \cdot)-\pi(\cdot)\| = 0.0636$, and
$\sup_{x\in {\cal X}} \|P^2(x, \cdot)-\pi(\cdot)\| = 0.00636$,
but $0.00636 > (0.0636)^2$.  On the other hand,
some authors instead define total variation distance as {\it twice} the
value used here, in which case the factor of~2 in
Proposition~\ref{P3}(e) is not written explicitly.
\end{remark}

\subsection{Asymptotic Convergence}
\label{sec-asymptotic}
\label{S3.2}
\label{S3.2}

Even if a Markov chain has stationary distribution $\pi(\cdot)$, it
may still fail to converge to stationarity:

\begin{ex}
Suppose ${\state} = \{1,2,3\}$, with $\pi\{1\}=\pi\{2\}=\pi\{3\}=1/3$.
Let $P(1,\{1\}) = P(1,\{2\}) = P(2,\{1\}) = P(2,\{2\}) = 1/2$,
and $P(3,\{3\}) = 1$.
Then $\pi(\cdot)$ is stationary.  However, if $X_0=1$, then
$X_n \in \{1,2\}$ for all
$n$, so $P(X_n=3)=0$ for all $n$, so $P(X_n=3) \not\to \pi\{3\}$, and
the distribution of $X_n$ does not converge to $\pi(\cdot)$.  (In fact,
here the stationary distribution is not {\it unique}, and the
distribution of $X_n$ converges to a different stationary distribution
defined by $\pi\{1\} = \pi\{2\} = 1/2$.)
\end{ex}

The above example is ``reducible'', in that the chain can never get
from state~1 to state~3, in any number of steps.  Now, the classical
notion of ``irreducibility'' is that the chain has positive probability
of eventually reaching
any state from any other state, but if $\X$ is uncountable then that
condition is impossible.  Instead, we demand the weaker condition of
$\phi$-irreducibility:\eject

\begin{defn}
A chain is {\it $\phi$-irreducible} if there exists a non-zero
$\sigma$-finite
measure $\phi$ on ${\state}$ such that for all $A \subseteq \state$
with $\phi(A) > 0$, and for all $x \in \state$,
there exists a positive integer $n=n(x,A)$ such that $P^{n}(x,A)>0$.
\end{defn}

For example, if $\phi(A) = \delta_{x_*}(A)$, then this requires that
$x_*$ has positive probability of eventually being reached from any
state $x$.  Thus, if a chain has any one state which is reachable from
anywhere (which on a finite state space is equivalent to being {\it
indecomposible}), then it is $\phi$-irreducible.  However, if $\X$
is uncountable then often $P(x,\{y\})=0$ for all $x$ and $y$.  In that
case, $\phi(\cdot)$ might instead be e.g.\ Lebesgue measure on $\IR^d$,
so that $\phi(\{x\})=0$ for all singleton sets, but such that all subsets
$A$ of positive Lebesgue measure are eventually reachable with positive
probability from any $x\in\X$.

\def\bx{{\bf x}}
\def\by{{\bf y}}
\def\b0{{\bf 0}}

\begin{Running Example}
Here we introduce a running example, to which we shall return several times.
Suppose that $\pi(\cdot)$ is a probability measure having
unnormalised density function $\pi_u$ with respect to
$d$-dimensional Lebesgue measure. Consider the Metropolis-Hastings
algorithm for $\pi_u$ with proposal density $q(\bx,\cdot )$ with respect to
$d$-dimensional Lebesgue measure. Then if $q(\cdot ,\cdot )$
is positive and continuous on $\IR^d \times \IR^d$,
and $\pi_u$ is finite everywhere, then the algorithm
is $\pi$-irreducible. Indeed, let $\pi(A)>0$.  Then
there exists $R >0$ such that $\pi(A_R)>0$,
where $A_R = A \cap B_R(\b0)$, and $B_R(\b0)$
represents the ball of radius $R$ centred at $\b0$.
Then by continuity, for any $\bx \in \IR^d$,
$\inf_{\by \in A_R} \min\{q(\bx, \by), q(\by, \bx)\}
\ge \epsilon$ for some $\epsilon > 0$, and thus we have
(assuming $\pi_u(\bx)>0$, otherwise $P(\bx,A)>0$ follows immediately) that
$$
P(\bx,A)  \ge
P(\bx,A_R)  \ge  \int_{A_R} q(\bx, \by) \, \min \left[
1,{ \pi_u(\by) \, q(\by, \bx) \over \pi_u(\bx) \, q(\bx, \by)} \right]
\, d\by
$$
$$
\ge \ \epsilon \
\hbox{Leb}\big(\{ \by \in A_R:\ \pi_u(\by) \ge \pi_u(\bx)\}\big)
\, + \,
{\epsilon \, K \over \pi_u(\bx )} \
\pi \big(\{ \by \in A_R:\ \pi_u(\by) < \pi_u(\bx)\} \big)
 ,
$$
where $K = \int_\X \pi_u(\bx) \, d\bx > 0$.
Since $\pi(\cdot)$ is absolutely continuous with respect to
Lebesgue measure, and since $Leb(A_R)>0$, it follows that
the terms in this final sum
cannot both be $0$, so that we must have $P(x,A) > 0$.
Hence, the chain is $\pi$-irreducible.
\end{Running Example}

Even $\phi$-irreducible chains might not converge in distribution,
due to periodicity problems, as in the following simple example.

\begin{ex}
\label{E2}
Suppose again ${\state} = \{1,2,3\}$, with $\pi\{1\}=\pi\{2\}=\pi\{3\}=1/3$.
Let $P(1,\{2\}) = P(2,\{3\}) = P(3,\{1\}) = 1$.
Then $\pi(\cdot)$ is stationary, and the chain is $\phi$-irreducible
[e.g.\ with $\phi(\cdot) = \delta_1(\cdot)$].  However, if $X_0=1$
(say), then $X_n = 1$ whenever
$n$ is a multiple of~3, so $P(X_n=1)$ oscillates between 0 and 1,
so again $P(X_n=1) \not\to \pi\{3\}$, and there is again no convergence
to $\pi(\cdot)$.
\end{ex}

To avoid this problem, we require {\it aperiodicity}, and we adopt the
following definition (which suffices for the $\phi$-irreducible chains
with stationary distributions that we shall study; for more general
relationships see e.g.\ Meyn and Tweedie~\cite{MT}, Theorem 5.4.4):

\begin{defn}
A Markov chain with stationary distribution $\pi(\cdot)$
is {\it aperiodic} if there do not exist $d \ge 2$ and
disjoint subsets ${\state}_1,{\state}_2,\ldots,{\state}_d \subseteq
\state$
with\break $P(x,{\state}_{i+1})=1$ for all $x\in{\state}_i$ ($1 \le i \le d-1$),
and $P(x,{\state}_{1})=1$ for all $x\in{\state}_d$, such that
$\pi({\state}_1)>0$ (and hence $\pi({\state}_i)>0$ for all $i$).
(Otherwise, the chain is {\it periodic}, with {\it period} $d$,
and {\it periodic decomposition} $\X_1,\ldots,\X_d$.)
\end{defn}

\begin{Running Example, Continued}
Here we return to the Running Example introduced above, and demonstrate
that no additional assumptions are necessary to ensure aperiodicity.
To see this, suppose that $\X_1$ and $\X _2$ are disjoint subsets of $\X$
both of positive $\pi $ measure, with $P(\bx , \X _2) = 1$
for all $ \bx \in
\X_1$. But just take any $\bx \in \X_1$, then  since $\X_1$ must
have positive Lebesgue measure,
$$
P(\bx ,\X_1)
 \ge  \int_{\by \in \bX_1} q(\bx ,\by ) \, \alpha(\bx,\by) \, d\by
 >  0
$$
for a contradiction. Therefore aperiodicity must hold.
(It is possible to demonstrate similar results for other MCMC algorithms,
such as the Gibbs sampler, see e.g.\ Tierney~\cite{tierney}.  Indeed, it
is rather rare for MCMC algorithms to be periodic.)
\end{Running Example, Continued}

Now we can state the main asymptotic convergence theorem, whose proof is
described in Section~\ref{S4}.  (This theorem assumes that the state space's
$\sigma$-algebra is {\it countably generated}, but this is a very weak
assumption which is true for e.g.\ any countable state space, or
any subset of $\IR^d$ with the usual
Borel $\sigma$-algebra, since that $\sigma$-algebra is generated by the
balls with rational centers and rational radii.)

\sc{thm}{3}
\begin{thm}
\label{asympconvthm}\label{T4}
If a Markov chain on a state space with countably generated $\sigma$-algebra
is $\phi$-irreducible and aperiodic, and has
a stationary distribution $\pi(\cdot)$, then for $\pi$-a.e.\ $x\in\state$,
$$
\lim_{n\to\infty} \| P^n(x,\cdot) - \pi(\cdot) \|  =  0.
$$
In particular, $\lim_{n\to\infty} P^n(x,A) = \pi(A)$ for
all measurable $A \subseteq \state$.
\end{thm}

\sc{Fact}{4}
\begin{Fact}
\label{SLLN}\label{F5}
In fact, under the conditions of Theorem~\ref{asympconvthm},
if $h:\state\to{\bf R}$ with $\pi(|h|) < \infty$, then
a ``strong law of large numbers'' also holds (see
e.g.\ Meyn and Tweedie~\cite{MT}, Theorem~17.0.1), as follows:
\be
\lim_{n \to \infty } (1/n) \, \sum _{i=1}^n h(X_i)  =  \pi(h)
\quad w.p.~1.
\label{SLLNeqn}\label{e6}
\ee
\end{Fact}

Theorem~\ref{asympconvthm} requires that the chain be
$\phi$-irreducible and aperiodic, and have stationary distribution
$\pi(\cdot)$.  Now, MCMC algorithms are created precisely so that
$\pi(\cdot)$ is stationary, so this requirement is not a problem.
Furthermore, it is usually straightforward to verify that chain is
$\phi$-irreducible, where e.g.\ $\phi$
is Lebesgue measure\vadjust{\eject} on an
appropriate region.  Also, aperiodicity almost always holds, e.g.\
for virtually any Metropolis algorithm or Gibbs sampler.  Hence,
Theorem~\ref{asympconvthm} is widely applicable to MCMC algorithms.

It is worth asking why the convergence in Theorem~\ref{asympconvthm}
is just from $\pi$-a.e.\ $x\in\X$.  The problem is that the chain may
have unpredictable behaviour on a ``null set'' of $\pi$-measure~0,
and fail to converge there.  Here is a simple example due to C.~Geyer
(personal communication):

\begin{ex}
\label{harrisex}%
Let $\X = \{1,2,\ldots\}$.  Let $P(1,\{1\}) = 1$, and for $x \ge 2$,
$P(x,\{1\}) = 1/x^2$ and $P(x,\{x+1\}) = 1-(1/x^2)$.
Then chain has stationary distribution $\pi(\cdot) = \delta_1(\cdot)$,
and it is $\pi$-irreducible and aperiodic.
On the other hand, if $X_0 = x \ge 2$, then $\P[X_n=x+n$ for all $n]
= \prod_{j=x}^\infty (1-(1/j^2)) > 0$, so that
$\|P^n(x,\cdot) - \pi(\cdot)\| \not\to 0$.
Here Theorem~\ref{asympconvthm} holds for $x=1$ which is indeed
$\pi$-a.e.\ $x\in\X$, but it does not hold for $x \ge 2$.
\end{ex}

\begin{remark}
The transient behaviour of the chain on the null set in
Example~\ref{harrisex} is not accidental.  If instead the chain
converged on the null set to some other stationary distribution, but
still had positive probability of escaping the null set (as it must to be
$\phi$-irreducible), then with probability~1 the chain would eventually
exit the null set, and would thus converge to $\pi(\cdot)$
from the null set after all.
\end{remark}

It is reasonable to ask under what circumstances the conclusions
of Theorem~\ref{asympconvthm} will hold for all $x\in\X$, not
just $\pi$-a.e.  Obviously, this will hold if the transition kernels
$P(x,\cdot)$ are all absolutely continuous with respect to $\pi(\cdot)$
(i.e., $P(x,dy) = p(x,y) \, \pi(dy)$ for some function $p:\X\times\X
\to [0,\infty)$), or for any Metropolis algorithm whose proposal
distributions $Q(x,\cdot)$ are absolutely continuous with respect to
$\pi(\cdot)$. It is also easy to see that this will hold for our Running
Example described above.  More generally, it suffices that the chain be
{\it Harris recurrent}, meaning that for all $B\subseteq\state$ with
$\pi(B)>0$, and all $x\in\state$, the chain will eventually reach $B$
from $x$ with probability~1, i.e.\ $\P[\exists \, n : \ X_n \in B \, |
\, X_0=x] = 1$.  This condition is stronger than $\pi$-irreducibility
(as evidenced by Example~\ref{harrisex}); for further discussions
of this see e.g.\ Orey~\cite{orey}, Tierney~\cite{tierney}, Chan and
Geyer~\cite{changeyer}, and \cite{harrismet}.

Finally, we note that periodic chains occasionally arise in MCMC (see
e.g.\ Neal~\cite{radfordnonrev}), and much of the theory can be applied
to this case.  For example, we have the following.

\sc{Corollary}{5}
\begin{Corollary}
\label{periodiccor}\label{C5}
If a Markov chain is $\phi$-irreducible, with period $d \ge 2$, and has
a stationary distribution $\pi(\cdot)$, then for $\pi$-a.e.\ $x\in\state$,
\be
\lim_{n\to\infty} \Big\| (1/d) \sum_{i=n}^{n+d-1} P^i(x,\cdot) - \pi(\cdot)
\Big\| = 0 ,
\label{periodicconveqn}\label{e7}
\ee
and also the strong law of large numbers~{\rm(\ref{SLLNeqn})} continues
to hold without change.
\end{Corollary}\vspace*{-24pt}\eject

\begin{proof}
Let the chain have periodic decomposition $\X_1,\ldots,\X_d \subseteq \X$,
and let $P'$ be the $d$-step chain $P^d$ restricted to the state
space $\X_1$.  Then $P'$ is $\phi$-irreducible and aperiodic on
$\X_1$, with stationary distribution $\pi'(\cdot)$ which satisfies
that $\pi(\cdot) = (1/d) \sum_{j=0}^{d-1} (\pi' \, P^j)(\cdot)$.
Now, from Proposition~\ref{tvprop}(c),
it suffices to prove the Corollary when $n=md$
with $m\to\infty$, and for simplicity we assume without loss of
generality that $x \in \X_1$.
From Proposition~\ref{tvprop}(d), we have
$\| P^{md+j}(x,\cdot) - (\pi' \, P^j)(\cdot) \|
\le \| P^{md}(x,\cdot) - \pi'(\cdot) \|$ for $j \in \IN$.
Then, by the triangle inequality,
$$
\Big\| (1/d) \sum_{i=md}^{md+d-1} P^i(x,\cdot) - \pi(\cdot) \Big\|
 =  \Big\| (1/d) \sum_{j=0}^{d-1} P^{md+j}(x,\cdot)
- (1/d) \sum_{j=0}^{d-1} (\pi' \, P^j)(\cdot) \Big\|
$$
$$
 \le  (1/d) \sum_{j=0}^{d-1} \| P^{md+j}(x,\cdot)
- (\pi' \, P^j)(\cdot) \|
 \le  (1/d) \sum_{j=0}^{d-1} \| P^{md}(x,\cdot)
- \pi'(\cdot) \|
 .
$$
But applying Theorem~\ref{asympconvthm} to $P'$, we obtain that
$\lim_{m\to\infty} \| P^{md}(x,\cdot) - \pi'(\cdot) \| = 0$ for
$\pi'$-a.e.\ $x \in \X_1$, thus giving the first result.

To establish~(\ref{SLLNeqn}), let $\Pbar$ be the transition kernel for
the Markov chain on $\X_1 \times \ldots \times \X_d$ corresponding
to the sequence
$\{(X_{md}, X_{md+1},\ldots,X_{md+d-1})\}_{m=0}^\infty$, and let
$\hbar(x_0,\ldots,x_{d-1}) = (1/d)(h(x_0)+\ldots+h(x_{d-1}))$.  Then
just like $P'$, we see that $\Pbar$ is $\phi$-irreducible and aperiodic,
with stationary distribution given by
$$
\pibar  =  \pi' \times (\pi' P) \times (\pi' P^2) \times
\ldots \times (\pi' P^{d-1})  .
$$
Applying Fact~\ref{SLLN} to $\Pbar$ and $\hbar$ establishes
that~(\ref{SLLNeqn}) holds without change.
\end{proof}

\begin{remark}
By similar methods, it follows that~(\ref{SLLN}) also remains true in
the periodic case, i.e.\ that
$$
\lim_{n \to \infty } (1/n) \, \sum _{i=1}^n h(X_i)  =  \pi(h)
\quad w.p.~1
$$
whenever $h:\state\to{\bf R}$ with $\pi(|h|) < \infty$,
provided the Markov chain is $\phi$-irreducible and countably generated,
without any assumption of aperiodicity.
In particular, both~(\ref{periodicconveqn}) and~(\ref{SLLN}) hold
(without further assumptions re periodicity)
for any irreducible (or indecomposible)
Markov chain on a {\it finite} state space.
\end{remark}

A related question for periodic chains, not considered here, is to
consider quantitative bounds on the difference of average distributions,
\[
\Bigg\| (1/n) \sum_{i=1}^{n} P^i(x,\cdot) - \pi(\cdot) \Bigg\|,
\]
through the
use of {\it shift-coupling}; see Aldous and Thorisson~\cite{aldthor},
and \cite{shift}.
\eject

\subsection{Uniform Ergodicity}
\label{S3.3}

Theorem~\ref{asympconvthm} implies asymptotic convergence to
stationarity, but does not say anything about the {\it rate} of this
convergence.  One ``qualitative'' convergence rate property is uniform
ergodicity:

\begin{defn}
A Markov chain having stationary distribution $\pi(\cdot)$
is {\it uniformly ergodic} if
$$
\|P^n(x,\cdot) - \pi(\cdot)\|
 \le  M \, \rho^n  ,
\qquad n=1,2,3,\ldots
$$
for some $\rho<1$ and $M<\infty$.
\end{defn}

One equivalence of uniform ergodicity is:

\sc{Proposition}{6}
\begin{Proposition}
\label{P7}
A Markov chain with stationary distribution
$\pi(\cdot)$ is uniformly ergodic if and only if
$\sup_{x\in\X} \|P^n(x,\cdot) - \pi(\cdot)\| < 1/2$ for some $n\in\IN$.
\end{Proposition}

\begin{proof}
If the chain is uniformly ergodic, then
$$
\lim_{n\to\infty} \, \sup_{x\in\X} \, \|P^n(x,\cdot) - \pi(\cdot)\|
 \le
\lim_{n\to\infty} M \, \rho^n
 =  0 ,
$$
so $\sup_{x\in\X} \|P^n(x,\cdot) - \pi(\cdot)\| < 1/2$ for all
sufficiently large $n$.  Conversely, if $\sup_{x\in\X} \|P^n(x,\cdot)
- \pi(\cdot)\| < 1/2$ for some $n\in\IN$, then in the notation of
Proposition~\ref{tvprop}(e), we have that $d(n) \equiv \beta < 1$,
so that for all $j\in\IN$, $d(jn) \le \big( d(n) \big)^j = \beta^j$.
Hence, from Proposition~\ref{tvprop}(c),
\beq
\|P^m(x,\cdot) - \pi(\cdot)\|& \le&
\|P^{\lfloor m/n \rfloor n}(x,\cdot) - \pi(\cdot)\| \le
\half d\left( \lfloor m/n \rfloor n \right)\\
&\le& \beta^{\lfloor m/n \rfloor}
\le \beta^{-1} \, \left(\beta^{1/n}\right)^m  ,
\eeq
so the chain is uniformly ergodic with $M = \beta^{-1}$ and
$\rho = \beta^{1/n}$.
\end{proof}

\begin{remark}
The above Proposition of course continues to hold if we replace $1/2$ by
$\delta$ for any $0<\delta<1/2$.  However, it is false for
$\delta \ge 1/2$.  For example, if $\X = \{1,2\}$, with $P(1,\{1\}) =
P(2,\{2\}) = 1$, and $\pi(\cdot)$ is uniform on $\X$, then
$\|P^n(x,\cdot) - \pi(\cdot)\| = 1/2$ for all $x\in\X$ and $n\in\IN$.
\end{remark}

To develop further conditions which ensure uniform ergodicity, we require
a definition.

\begin{defn}
A subset $C \subseteq \state$ is {\it small} (or,
$(n_0,\epsilon,\nu)$-small) if there exists a positive
integer $n_0$, $\epsilon >0$, and a probability measure $\nu(\cdot)$ on
$\X$ such that the following {\it minorisation condition} holds:
\be
P^{n_0}(x,\cdot)  \ge  \epsilon \, \nu(\cdot)
\qquad x\in C ,
\label{minorcond}
\label{e8}
\ee
i.e.\ $P^{n_0}(x,A) \ge \epsilon \, \nu(A)$ for all $x\in C$ and
all measurable $A \subseteq \X$.
\end{defn}

\begin{remark}
Some authors (e.g.\ Meyn and Tweedie~\cite{MT}) also
require that $C$ have positive
stationary measure, but for simplicity we don't explicitly require that
here.  In any case, $\pi(C)>0$ follows under the additional
assumption of the drift
condition~(\ref{unidriftcond}) considered in the
next section.
\end{remark}

Intuitively, this condition means that all of the $n_0$-step
transitions from within $C$, all have an ``$\epsilon$-overlap'', i.e.\
a component of size $\epsilon$ in common.  (This concept goes back to
Doeblin~\cite{doeblin}; for further background, see e.g.\ \cite{doob},
\cite{athreya}, \cite{nummelin}, \cite{asmussen}, and~\cite{MT}; for
applications to convergence rates see e.g.\ \cite{MT2}, \cite{ros95},
\cite{jamesstein}, \cite{bruns}, \cite{RT}, \cite{douc}, \cite{ros02}.)
We note that if $\X$ is countable, and if
\be
\epsilon_{n_0}  \equiv  \sum_{y\in\X} \inf_{x\in C} P^{n_0}(x,\{y\})
 >  0  ,
\label{discepsilon}\label{e9}
\ee
then $C$ is $(n_0,\epsilon_{n_0},\nu)$-small where
$\nu\{y\} = \epsilon_{n_0}^{-1} \inf_{x\in C} P^{n_0}(x,\{y\})$.
(Furthermore, for an irreducible (or
just indecomposible) and aperiodic chain on a {\it finite} state space,
we always have $\epsilon_{n_0}>0$ for sufficiently large $n_0$ (see
e.g.~\cite{eigen}), so this method always applies in principle.)
Similarly, if the transition probabilities
have densities with respect to some measure $\eta(\cdot)$, i.e.\ if
$P^{n_0}(x,dy) = p_{n_0}(x,y) \, \eta(dy)$, then we can take
$\epsilon_{n_0} = \int_{y\in\X} \big( \inf_{x\in\X} p_{n_0}(x,y) \big)
\, \eta(dy)$.

\begin{remark}
As observed in~\cite{pseudosmall}, small-set conditions of
the form $P(x,\cdot) \ge \epsilon \, \nu(\cdot)$ for all $x\in C$, can be
replaced by {\it pseudo-small} conditions of the form $P(x,\cdot) \ge
\epsilon \, \nu_{xy}(\cdot)$ and $P(y,\cdot) \ge \epsilon \, \nu_{xy}(\cdot)$
for all $x,y\in C$, without affecting any bounds which use pairwise coupling
(which includes all of the bounds considered here before
Section~\ref{S5}.  Thus, all of the
results stated in this section remain true without change if ``small set''
is replaced by ``pseudo-small set'' in the hypotheses.  For ease of
exposition, we do not emphasise this point herein.
\end{remark}

The main result guaranteeing uniform ergodicity, which goes back to
Doeblin~\cite{doeblin} and Doob~\cite{doob} and in some sense even to
Markov~\cite{markovref}, is the following.

\sc{thm}{7}
\begin{thm}
\label{unifthm}
Consider a Markov chain with invariant probability distribution
$\pi(\cdot)$.  Suppose the minorisation
condition~{\rm(\ref{minorcond})}
is satisfied for some $n_0\in\IN$ and $\epsilon>0$ and probability
measure $\nu(\cdot)$, in the special case $C=\X$ (i.e., the entire
state space is small).  Then the chain is uniformly ergodic, and in fact
$\|P^n(x,\cdot)-\pi(\cdot)\| \le (1-\epsilon)^{\lfloor n/n_0 \rfloor}$
for all $x\in\X$, where $\lfloor r \rfloor$ is the greatest integer not
exceeding~$r$.
\end{thm}

Theorem~\ref{unifthm} is proved in Section~\ref{S4}.  We note also that
Theorem~\ref{unifthm} provides a {\it quantitative} bound on
the distance to stationarity $\|P^n(x,\cdot)-\pi(\cdot)\|$, namely
that it must be $\le (1-\epsilon)^{\lfloor n/n_0 \rfloor}$.  Thus,
once $n_0$ and $\epsilon$ are known, we can find $n_*$ such that, say,
$\|P^{n_*}(x,\cdot)-\pi(\cdot)\| \le 0.01$, a fact which can be applied
in certain MCMC contexts (see e.g.~\cite{ros93}).  We can then say that
$n_*$ iterations ``suffices for convergence'' of the Markov chain.
On a discrete state space, we have that $\|P^n(x,\cdot)-\pi(\cdot)\|
\le (1-\epsilon_{n_0})^{\lfloor n/n_0 \rfloor}$ with $\epsilon_{n_0}$
as in~(\ref{discepsilon}).\eject

\begin{Running Example, Continued}
Recall our Running Example, introduced above.
Since we have imposed strong continuity
conditions on $q$, it is natural to conjecture that compact sets are small.
However this is not true without extra regularity conditions. For instance,
consider dimension $d=1$, and
suppose that $\pi_u(x) = \one_{0<|x|<1} |x|^{-1/2}$, and let $q(x,y) \propto
\exp\{-(x-y)^2/2\}$, then it is easy to check that any neighbourhood of $0$ is
not small. However in the general setup of our Running Example,
all compact sets on which $\pi_u$ is bounded are
small. To see this, suppose $C$ is a compact set on which $\pi_u$
is bounded by $k<\infty $. Let $\bx \in C$, and let $D$ be any compact set of
positive Lebesgue and $\pi $ measure, such that
$\inf_{\bx \in C, \by \in D}q(\bx ,\by )=\epsilon >0$
for all $\by \in D$.  We then have,
$$
P(\bx , d \by)  \ge  q(\bx , \by ) \, d\by \, \min \left\{1, {\pi_u(\by )
\over \pi_u(\bx)}\right\}
 \ge  \epsilon \, d\by \min \left\{1, {\pi_u(\by ) \over k} \right\}
 ,
$$
which is a positive measure independent of $\bx$.
Hence, $C$ is small.
(This example also shows that if $\pi_u$ is continuous, the state
space $\X$ is compact, and $q$ is continuous and positive, then $\X$
is small, and so the chain must be uniformly ergodic.)
\end{Running Example, Continued}

If a Markov chain is {\it not} uniformly ergodic (as few MCMC algorithms
on unbounded state spaces are), then Theorem~\ref{unifthm} cannot
be applied.  However, it is still of great importance, given a Markov
chain kernel $P$ and an initial state $x$, to be able to
find $n_*$ so that, say, $\|P^{n_*}(x,\cdot)-\pi(\cdot)\| \le 0.01$.
This issue is discussed further below.

\subsection{Geometric ergodicity}
\label{sec-geometric}
\label{S3.4}

A weaker condition than uniform ergodicity is geometric ergodicity, as
follows (for background and history, see e.g.\ Nummelin~\cite{nummelin},
and Meyn and Tweedie~\cite{MT}):

\begin{defn}
A Markov chain with stationary distribution $\pi(\cdot)$
is {\it geometrically ergodic} if
$$
\|P^n(x,\cdot) - \pi(\cdot)\|
 \le  M(x) \, \rho^n  ,
\qquad n=1,2,3,\ldots
$$
for some $\rho<1$, where $M(x)<\infty$ for $\pi$-a.e.\ $x\in\state$.
\end{defn}
The difference between geometric ergodicity and uniform ergodicity is
that now the constant $M$ may depend on the initial state~$x$.

Of course, if the state space $\X$ is {\it finite}, then all
irreducible and aperiodic Markov chains are geometrically (in fact,
uniformly) ergodic.  However, for infinite $\X$ this is not the case.
For example, it is shown by Mengersen and Tweedie~\cite{mengersen} (see
also~\cite{RTmet}) that a symmetric random-walk Metropolis algorithm
is geometrically ergodic essentially if and only if $\pi(\cdot)$ has
finite exponential moments.  (For chains which are not geometrically
ergodic, it is possible also to study {\it polynomial ergodicity}, not
considered here; see Fort and Moulines~\cite{fortmoulines}, and Jarner
and Roberts~\cite{jarner}.)  Hence, we now discuss conditions which
ensure geometric ergodicity.

\begin{defn}
Given Markov chain transition probabilities $P$ on a state space $\X$,
and a measurable
function $f:\X\to\IR$, define the function $Pf : \X \to \IR$ such
that $(Pf)(x)$ is the conditional expected value of $f(X_{n+1})$, given
that $X_n = x$.  In symbols, $(Pf)(x) = \int_{y\in\X} f(y) \, P(x,dy)$.
\end{defn}

\begin{defn}
A Markov chain satisfies a {\it drift condition} (or, univariate
geometric drift condition) if
there are constants $0<\lambda <1$ and $b<\infty$,
and a function $V : {\state} \to [1,\infty]$, such that
\be
PV  \le  \lambda V + b \one_C,
\label{unidriftcond}\label{e10}
\ee
i.e.\ such that $\int_{\state } P(x, dy ) V(y) \le \lambda V(x) + b{\bf
1}_C (x)$ for all $x\in\X$.
\end{defn}

The main result guaranteeing geometric ergodicity is the following.


\sc{thm}{8}
\begin{thm}
\label{geomthm}\label{T9}
Consider a $\phi$-irreducible, aperiodic Markov chain with
stationary distribution~$\pi(\cdot)$.  Suppose the minorisation
condition~{\rm(\ref{minorcond})} is
satisfied for some $C \subset \X$ and $\epsilon>0$ and probability measure
$\nu(\cdot)$.  Suppose further that the drift
condition~{\rm(\ref{unidriftcond})}
is satisfied for some constants $0<\lambda <1$ and $b<\infty$ , and a
function $V : {\state} \to [1,\infty]$
with $V(x)<\infty$ for at least one (and hence for $\pi$-a.e.) $x\in\X$.
Then then chain is geometrically ergodic.
\end{thm}

Theorem~\ref{geomthm} is usually proved by complicated analytic
arguments (see e.g.\ \cite{nummelin}, \cite{MT}, \cite{doss}).
In Section~\ref{S4}, we describe a proof of Theorem~\ref{geomthm}
which uses direct coupling constructions instead.  Note also that
Theorem~\ref{geomthm} provides no {\it quantitative} bounds on $M(x)$
or $\rho$, though this is remedied in
Theorem~\ref{quantbound} below.

\sc{Fact}{9}
\begin{Fact}
\label{Vunifremark}\label{F10}
In fact, it follows from Theorems~15.0.1,~16.0.1, and~14.3.7 of Meyn
and Tweedie~\cite{MT}, and Proposition~1 of~\cite{hybrid}, that
the minorisation condition~(\ref{minorcond}) and
drift condition~(\ref{unidriftcond})
of Theorem~\ref{geomthm} are equivalent
(assuming $\phi$-irreducibility and aperiodicity)
to the apparently stronger property of
``$V$-uniform ergodicity'', i.e.\ that there is $C<\infty$
and $\rho<1$ such that
$$
\sup_{|f| \le V} |P^nf(x) - \pi(f)|  \le  C \, V(x) \, \rho^n,
\qquad x \in \X ,
$$
where $\pi(f) = \int_{x\in\X} f(x) \, \pi(dx)$.  That is,
we can take $\sup_{|f| \le V}$ instead of just $\sup_{0<f<1}$
(compare Proposition~\ref{tvprop} parts~(a) and~(b)), and we can
let $M(x) = C \, V(x)$ in the geometric ergodicity bound.
Furthermore, we always have $\pi(V)<\infty$.  (The term
``$V$-uniform ergodicity'', as used in~\cite{MT}, perhaps also implies that
$V(x)<\infty$ for all $x\in\X$, rather than just for $\pi$-a.e.\ $x\in\X$,
though we do not consider that distinction further here.)
\end{Fact}\vspace*{-12pt}\eject

\begin{opennew}
Can direct coupling methods, similar to those used below to prove
Theorem~{\rm\ref{geomthm}}, also be used to provide an alternative proof
of Fact~{\rm\ref{Vunifremark}}?
\end{opennew}

\begin{ex}
Here we consider a simple example of
geometric ergodicity of Metropolis algorithms on $\IR$ (see Mengersen
and Tweedie~\cite{mengersen}, and~\cite{RTmet}).  Suppose that
$\X = \IR ^+$ and $\pi_u(x) = e^{-x}$. We will use a symmetric (about $x$)
proposal
distribution $q(x, y)= q(|y-x|)$ with support contained in $[x-a, x+a]$.
In this simple situation,
a natural drift function to take is $V(x) = e^{cx}$ for some
$c>0$.  For $x \ge a$, we compute:
\beq
PV(x) & = &
\int_{x-a}^x V(y)q(x,y)dy + \int_x^{x+a}V(y)q(x,y) dy
{\pi_u(y) \over \pi_u(x)}\\
&&{}+ V(x) \int_x^{x+a} q(x,y)dy (1 - \pi_u(y) /\pi_u(x)) .
\eeq
By the symmetry of $q$, this can be written as
$$
\int_x^{x+a}
I(x,y) \, q(x, y) \, dy ,
$$
where
\beq
I(x,y) &=& {V(y) \pi_u(y) \over \pi_u(x)} + V(2x-y) + V(x) \left( 1-
{\pi_u(y) \over \pi_u(x)}\right)\\
&=& e^{cx} \left[
e^{(c-1)u} + e^{-cu} + 1 - e^{-u}
\right]
= e^{cx}\left[
2 - (1+ e^{(c-1)u})(1-e^{-cu})
\right],
\eeq
and where $u = y - x$.
For $c<1$, this is equal to $2(1-\epsilon ) V(x)$
for some positive constant $\epsilon$.
Thus in this case we have shown that for all $x>a$
$$
PV(x)  \le  \int_x^{x+a} 2V(x) (1-\epsilon ) q(x,y) dy
 =  (1-\epsilon ) V(x).
$$
Furthermore, it is easy to show that $PV(x)$ is bounded on $[0,a]$ and that
$[0,a]$ is in fact a small set.  Thus,
we have demonstrated that the drift condition~(\ref{unidriftcond})
holds.  Hence, the algorithm is
geometrically ergodic by Theorem~\ref{geomthm}.
(It turns out that for such Metropolis algorithms,
a certain condition, which essentially requires
an exponential bound on the tail probabilities of $\pi(\cdot)$, is in fact
{\it necessary} for geometric ergodicity; see~\cite{RTmet}.)
\end{ex}

Implications of geometric ergodicity for central limit theorems are
discussed in Section~\ref{S5}.  In general, it believed by practitioners of
MCMC that geometric ergodicity is a useful property.  But does geometric
ergodicity really matter?  Consider the following examples.\eject

\begin{ex}
(\cite{bruns})
\label{indepexample}\label{E5}
Consider an independence sampler, with
$\pi(\cdot)$ an Ex\-po\-nen\-tial$(1)$ distribution, and $Q(x,\cdot)$ an
Exponential$(\lambda)$ distribution.  Then if $0 < \lambda \le 1$,
the sampler is geometrically ergodic, has central limit theorems (see
Section~\ref{S5}), and generally behaves fairly well even for very small
$\lambda$.  On the other hand, for $\lambda > 1$ the sampler fails
to be geometrically ergodic, and indeed for $\lambda \ge 2$ it fails
to have central limit theorems, and generally behaves quite poorly.
For example, the simulations in~\cite{bruns} indicate that
with $\lambda=5$, when started in stationarity and
averaged over the first million iterations,
the sampler will usually return an average
value of about $0.8$ instead of $1$, and then occasionally return a very
large value instead, leading to very unstable behaviour.  Thus, this
is an example where the property of geometric ergodicity does indeed
correspond to stable, useful convergence behaviour.
\end{ex}

However, geometric ergodicity does not always guarantee a useful Markov
chain algorithm, as the following two examples show.

\begin{ex}
(``Witch's Hat'', e.g.\ Matthews~\cite{matthews}) Let ${\state}
= [0,1]$, let $\delta = 10^{-100}$ (say), let $0<a<1-\delta$, and let
$\pi_u(\bx) = \delta + \one_{[a,a+\delta]}(\bx)$.  Then $\pi([a,a+\delta])
\approx 1/2$.  Now, consider running a typical Metropolis algorithm on
$\pi_u$.  Unless $X_0 \in [a,a+\delta]$, or the sampler gets ``lucky''
and achieves $X_n \in [a,a+\delta]$ for some moderate $n$, then the
algorithm will likely \un{miss} the tiny interval $[a,a+\delta]$ entirely,
over any feasible time period.  The algorithm will thus ``appear''
(to the naked eye or to any statistical test) to converge to the ${\rm
Uniform}({\state})$ distribution, even though ${\rm Uniform}({\state})$
is very different from $\pi(\cdot)$.  Nevertheless, this algorithm is
still geometrically ergodic (in fact uniformly ergodic).  So in this
example, geometric ergodicity does not guarantee a well-behaved sampler.
\end{ex}

\begin{ex}
Let $\state = {\bf R}$, and let $\pi_u(x) = 1/(1+x^2)$ be the
(unnormalised) density of the Cauchy distribution.
Then a random-walk Metropolis algorithm
for $\pi_u$ (with, say, $X_0=0$ and $Q(x,\cdot) = {\rm
Uniform}[x-1,x+1]$) is ergodic but is {\it not} geometrically ergodic.
And, indeed, this sampler has very slow, poor convergence properties.
On the other hand, let $\pi'_u(x) = \pi_u(x) \, \one_{|x| \le 10^{100}}$,
i.e.\ $\pi'_u$ corresponds to $\pi_u$ truncated at $\pm$ one googol.
Then the same random-walk Metropolis algorithm for
$\pi'_u$ {\it is} geometrically ergodic, in fact uniformly ergodic.
However, the two algorithms are {\it indistinguishable} when
run for any remotely feasible number of iterations.
Thus, this is an example where geometric ergodicity does not in any way
indicate improved performance of the algorithm.
\end{ex}

In addition to the above two examples, there are also numerous examples
of important Markov chains on {\it finite} state spaces
(such as the single-site Gibbs sampler for the Ising model
at low temperature
on a large but finite grid) which are irreducible and aperiodic,
and hence uniformly (and thus also geometrically) ergodic, but which
converge to stationarity extremely slowly.

The above examples illustrate a limitation of {\it qualitative}
convergence properties such as geometric ergodicity.  It is thus desirable
where possible to instead obtain {\it quantitative} bounds on Markov
chain convergence.  We consider this issue next.

\subsection{Quantitative Convergence Rates}
\label{S3.5}

In light of the above, we ideally want {\it quantitative} bounds on
convergence rates, i.e.\ bounds of the form $\|P^n(x,\cdot)-\pi(\cdot)\|
\le g(x,n)$ for some {\it explicit} function $g(x,n)$, which (hopefully)
is small for large $n$.  Such questions now have a substantial
history in MCMC, see e.g.\ \cite{MT2}, \cite{ros95},
\cite{jamesstein}, \cite{lund}, \cite{cowles1}, \cite{RT}, \cite{hobert1},
\cite{hobert2}, \cite{douc}, \cite{berenhaut}, \cite{ros02}, \cite{fort},
\cite{baxendale}, \cite{period}, \cite{periodnon}.

We here present a result from~\cite{ros02}, which follows as a special
case of~\cite{douc}; it is based on the approach of~\cite{ros95} while
also taking into account a small improvement from~\cite{RT}.

Our result requires a {\it bivariate drift condition} of the form
\be
\overline{P} h(x,y)
 \le  h(x,y)  /  \alpha  ,
\qquad (x,y) \notin C \times C
\label{bivardriftcond}
\label{e11}
\ee
for some function $h:\X\times\X\to[1,\infty)$ and some $\alpha>1$, where
$$
\overline{P} h(x,y)  \equiv
\int_\X \int_\X h(z,w) \, P(x,dz) \, P(y,dw)  .
$$
(Thus, $\overline{P}$ represents running two independent copies of the chain.)
Of course,~(\ref{bivardriftcond}) is closely related
to~(\ref{unidriftcond});
for example we have the following (see also~\cite{ros95}, and
Proposition~2 of~\cite{periodnon}):

\sc{Proposition}{10}
\begin{Proposition}
\label{unibiprop}\label{P10}
Suppose the univariate drift
condition~{\rm(\ref{unidriftcond})} is satisfied
for some $V:\X\to[1,\infty]$, $C \subseteq \X$, $\lambda<1$, and $b<\infty$.
Let $d = \inf_{x\in C^c} V(x)$.  Then if $d > [b/(1-\lambda)]-1$, then the
bivariate drift
condition~{\rm(\ref{bivardriftcond})} is satisfied for the
same $C$, with $h(x,y) = \half[V(x)+V(y)]$ and $\alpha^{-1} = \lambda +
b/(d+1) < 1$.
\end{Proposition}

\begin{proof}
If $(x,y) \notin C \times C$, then either $x \notin C$ or $y\notin C$
(or both), so $h(x,y) \ge (1 + d)/2$, and $PV(x)+PV(y) \le \lambda
V(x) + \lambda V(y) + b$.  Then
\beq
\overline{P} h(x,y) &=& \half[PV(x) + PV(y)]
\le \half[\lambda V(x) + \lambda V(y) + b]\\
&=& \lambda \, h(x,y) + b/2
\le \lambda \, h(x,y) + (b/2)[h(x,y)/((1+d)/2)]\\
&=& [\lambda + b/(1+d)] \, h(x,y) .
\eeq
Furthermore, $d > [b/(1-\lambda)]-1$ implies that $\lambda + b/(1+d) < 1$.
\end{proof}

Finally, we let
\be
B_{n_0}  =  \max\left[1, \
\alpha^{n_0} (1-\epsilon) \sup_{C \times C}
\overline{R} h\right],
\label{Bdef}\label{e12}
\ee
where for $(x,y) \in C \times C$,
$$
\overline{R}h(x,y)  =
\int_\X \int_\X (1-\epsilon)^{-2} h(z,w)
\, (P^{n_0}(x,dz) - \epsilon \nu(dz))
\, (P^{n_0}(y,dw) - \epsilon \nu(dw)).
$$

In terms of these assumptions, we state our result as follows.

\sc{thm}{11}
\begin{thm}
\label{quantbound}\label{T12}
Consider a Markov chain on a state space $\X$, having transition kernel
$P$.  Suppose there is $C \subseteq \X$, $h:\X\times\X \to [1,\infty)$,
a probability distribution $\nu(\cdot)$ on $\X$,
$\alpha>1$, $n_0\in\IN$, and $\epsilon>0$, such
that~{\rm(\ref{minorcond})}
and~{\rm(\ref{bivardriftcond})} hold.
Define $B_{n_0}$ by~{\rm(\ref{Bdef})}.
Then for any joint initial distribution $\L(X_0,X'_0)$, and
any integers $1 \le j \le k$, if $\{X_n\}$ and $\{X'_n\}$ are two copies
of the Markov chain started in the joint
initial distribution $\L(X_0,X'_0)$, then
$$
\| \L(X_k) - \L(X'_k) \|_{TV}
 \le (1-\epsilon)^j +
\alpha^{-k} \, (B_{n_0})^{j-1} \, \E[h(X_0,X'_0)] .
$$
In particular, by choosing $j = \lfloor rk \rfloor$ for sufficiently
small $r>0$, we obtain an explicit, quantitative convergence bound which
goes to~0 exponentially quickly as $k\to\infty$.
\end{thm}

Theorem~\ref{quantbound} is proved in Section~\ref{S4}.  Versions of
this theorem have been applied to various realistic MCMC algorithms,
including for versions of the variance components model described
earlier, resulting in bounds like $\|P^n(x,\cdot)-\pi(\cdot)\| < 0.01$
for $n=140$ or $n=3415$; see e.g.\ \cite{jamesstein}, and Jones and
Hobert~\cite{hobert2}.  Thus, while it is admittedly hard work to apply
Theorem~\ref{quantbound} to realistic MCMC algorithms, it is indeed
possible and often can establish rigorously that perfectly feasible
numbers of iterations are sufficient to ensure convergence.

\begin{remark}
For complicated Markov chains, it might be difficult to apply
Theorem~\ref{quantbound} successfully.  In such cases, MCMC
practitioners instead use ``convergence diagnostics'', i.e.\ do
statistical analysis of the realised output $X_1,X_2,\ldots$, to
see if the distributions of $X_n$ appear to be ``stable'' for large
enough $n$.  Many such diagnostics involve running the Markov chain
repeatedly from different initial states, and checking if the chains all
converge to approximately the same distribution (see e.g.\ Gelman and
Rubin~\cite{gelmanrubin}, and Cowles and Carlin~\cite{cowlescarlin}).
This technique often works well in practice.  However, it provides
no rigorous guarantees and can sometimes be fooled into prematurely
claiming convergence (see e.g.~\cite{matthews}), as is likely to happen
for the examples at the end of
Section~\ref{S3}.  Furthermore,
convergence diagnostics can also introduce bias into the resulting
estimates (see~\cite{biases}).  Overall, despite the extensive theory
surveyed herein, the ``convergence time problem'' remains largely
unresolved for practical application of MCMC.  (This is also the
motivation for ``perfect MCMC'' algorithms, originally developed by
Propp and Wilson~\cite{proppwilson} and not discussed here; for further
discussion see e.g.~Kendall and M\o ller~\cite{kendallmoller},
Th\"onnes~\cite{thonnes}, and Fill et al.~\cite{fmmr}.)
\end{remark}

\section{Convergence Proofs using Coupling Constructions}
\label{S4}

In this section, we prove some of the theorems stated earlier.  There are
of course many methods available for bounding convergence of Markov
chains, appropriate to various settings (see e.g.\ \cite{aldous},
\cite{diaconis}, \cite{sinclair}, \cite{sst}, \cite{stein}, and
Subsection~\ref{S5.4} herein), including the setting of large but
finite state spaces that often arises in computer science (see e.g.\
Sinclair~\cite{sinclair} and Randall~\cite{randall}) but is not our
emphasis here.  In this section, we focus on the method of {\it coupling},
which seems particularly well-suited to analysing MCMC algorithms on
general (uncountable) state spaces.  It is also particularly well-suited
to incorporating small sets (though small sets can also be combined with
{\it regeneration theory}, see e.g.\ \cite{athreya}, \cite{asmussen},
\cite{mykland}, \cite{hobertregen}).  Some of the proofs below are new,
and avoid many of the long analytic arguments of some previous proofs
(e.g.\ Nummelin~\cite{nummelin}, and Meyn and Tweedie~\cite{MT}).


\subsection{The Coupling Inequality}
\label{S4.1}

The basic idea of coupling
is the following.  Suppose we have two \it random
variables \rm $X$ and $Y$, defined jointly on some space $\X$.  If we
write $\L(X)$ and $\L(Y)$ for their respective probability
distributions, then we can write
\begin{eqnarray*}
\|\L(X)-\L(Y)\|  &= & \supl_A |P(X \in A) - P(Y \in A)| \\
&=& \supl_A |P(X \in A, X=Y) + P(X \in A, X \not= Y) \\
&&{} - P(Y \in A, Y=X) - P(Y \in A, Y \not =X) | \\
&=& \supl_A |P(X \in A, X \not= Y) - P(Y \in A, Y \not =X) | \\
&\le& P(X \not= Y)  ,
\end{eqnarray*}
so that
\be
\|\L(X)-\L(Y)\| \ \le P(X \not= Y)  .
\label{couplinginequality}
\label{e13}
\ee
That is, {\sl the variation distance between the laws of two
random variables is bounded by the probability that they are
unequal}.  For background, see e.g.\ Pitman~\cite{pitman},
Lindvall~\cite{lindvall}, and Thorisson~\cite{thorisson}.

\subsection{Small Sets and Coupling}
\label{sec-coupling}
\label{S4.2}

Suppose now that $C$ is a small set.  We shall use the following
coupling construction, which is essentially the ``splitting technique''
of Nummelin~\cite{nummelinpaper} and Athreya and Ney~\cite{athreya};
see also Nummelin~\cite{nummelin}, and Meyn and Tweedie~\cite{MT}.
The idea is to run two copies $\{X_n\}$ and $\{X'_n\}$ of the Markov
chain, each of which marginally follows the updating rules $P(x,\cdot)$,
but whose joint construction (using $C$) gives them as high a
probability as possible of becoming equal to each other.

\bigskip\noindent THE COUPLING CONSTRUCTION:

Start with $X_0=x$ and $X'_0 \sim \pi(\cdot)$, and $n=0$, and repeat the
following loop forever.

{\bf Beginning of Loop.} Given $X_n$ and $X'_n$:

\begin{enumerate}
\item[1.]
If $X_n=X'_n$, choose $X_{n+1} = X'_{n+1} \sim P(X_n,\cdot)$, and replace
$n$ by $n+1$.

\item[2.]
Else, if $(X_n,X'_n) \in C \times C$, then:

\hskip 2cm
{(a)} w.p.\ $\epsilon$, choose $X_{n+n_0}=X'_{n+n_0} \sim \nu(\cdot)$;

\hskip 2cm
{(b)} else, w.p.\ $1-\epsilon$, conditionally independently choose
$$
X_{n+n_0} \sim {1 \over 1-\epsilon} [ P^{n_0}(X_n,\cdot)
- \epsilon \, \nu(\cdot)] ,
$$
$$
X'_{n+n_0} \sim {1 \over 1-\epsilon}
[ P^{n_0}(X'_n,\cdot) - \epsilon \, \nu(\cdot)]
 .
$$

In the case $n_0>1$, for completeness go back and construct
$X_{n+1},\ldots,X_{n+n_0-1}$ from their correct conditional
distributions given $X_n$ and $X_{n+n_0}$, and similarly (and
conditionally independently) construct
$X'_{n+1},\ldots,X'_{n+n_0-1}$ from their correct conditional
distributions given $X'_n$ and $X'_{n+n_0}$.
In any case, replace $n$ by $n+n_0$.

\item[3.]
Else, conditionally independently choose
$X_{n+1} \sim P(X_n,\cdot)$ and
$X'_{n+1} \sim P(X'_n,\cdot)$, and replace $n$ by $n+1$.
\end{enumerate}

{\bf Then return to Beginning of Loop.}

\bigskip

Under this construction, it is easily checked that $X_n$ and $X'_n$
are each marginally updated according to the correct transition kernel $P$.
It follows that $\P[X_n\in A] = P^n(x,\cdot)$ and $\P[X'_n \in A] =
\pi(A)$ for all $n$. Moreover the two chains are run
independently until they both enter $C$ at which time the
minorisation {\it splitting} construction (step 2) is utilised.
Without such a construction, on uncountable state spaces, we
would not be able to ensure successful coupling of the two
processes.

The {\it coupling inequality} then says that $\|P^n(x,\cdot) - \pi(\cdot)\|
\le \P[X_n\not=X'_n]$.  The question is, can we use this to obtain useful
bounds on $\|P^n(x,\cdot) - \pi(\cdot)\|$?  In fact, we shall now provide
proofs (nearly self-contained) of all of the theorems stated earlier,
in terms of this coupling construction.  This allows for intuitive
understanding of the theorems, while also avoiding various analytic
technicalities of the previous proofs of some of these theorems.

\subsection{Proof of Theorem~\ref{unifthm}}
\label{S4.3}

In this case, $C=\X$, so every $n_0$ iterations we have probability
at least $\epsilon$ of making $X_n$ and $X'_n$ equal.  It follows that
if $n = n_0 m$, then $\P[X_n \not= X'_n] \le (1-\epsilon)^m$.  Hence,
from the coupling inequality,
$\|P^n(x,\cdot) - \pi(\cdot)\| \le (1-\epsilon)^m = (1-\epsilon)^{n/n_0}$
in this case.
It then follows from Proposition~\ref{tvprop}(c) that
$\|P^n(x,\cdot) - \pi(\cdot)\| \le (1-\epsilon)^{ \lfloor n/n_0 \rfloor}$
for any $n$.\qed

\subsection{Proof of Theorem~\ref{quantbound}}
\label{S4.4}

We follow the general outline of~\cite{ros02}.  We again begin by assuming
that $n_0=1$ in the minorisation condition for the small set $C$ (and
thus write $B_{n_0}$ as $B$), and indicate at the end what changes are
required if $n_0>1$.

Let
$$
N_k  =  \#\{m \ : \ 0 \le m \le k, \ (X_m,X'_m) \in C \times C \}  ,
$$
and let $\tau_1,\tau_2,\ldots$ be the times of the
successive visits of $\{(X_n,X'_n)\}$ to $C \times C$.
Then for any integer $j$ with $1 \le j \le k$,
\be
\P[X_k\not=X'_k]
 =
\P[X_k\not=X'_k, \ N_{k-1} \ge j]  +  \P[X_k\not=X'_k, \ N_{k-1} < j]
 .
\label{twoterms}\label{e14}
\ee

Now, the event $\{X_k\not=X'_k, \ N_{k-1} \ge j\}$ is contained in
the event that the first $j$ coin flips all came up tails.  Hence,
$\P[X_k\not=X'_k, \ N_{k-1} \ge j] \le (1-\epsilon)^j$.  which bounds
the first term in~(\ref{twoterms}).

To bound the second term in~(\ref{twoterms}), let
$$
M_k  =  \alpha^{k} B^{-N_{k-1}} h(X_k,X'_k) \, \one(X_k\not=X'_k) ,
\qquad k=0,1,2,\ldots
$$
(where $N_{-1}=0$).

\sc{Lemma}{12}
\begin{Lemma}
We have
$$
\E[M_{k+1} \mid X_0,\ldots,X_k,X'_0,\ldots,X'_k]  \le  M_k  ,
$$
i.e.\ $\{M_k\}$ is a {\it supermartingale}.
\end{Lemma}

\begin{proof}
If $(X_k, X'_k)\notin C\times C$, then $N_k=N_{k-1}$, so
\begin{eqnarray*}
&&\E[M_{k+1} \mid X_0,\ldots, X_k, X'_0,\ldots, X'_k]
\qquad\qquad\qquad \\
&&\qquad= \ \alpha^{k+1} B^{-N_{k-1}} \E[h(X_{k+1}, X'_{k+1})
\, {\bf 1}(X_{k+1}\neq X'_{k+1}) \mid
X_k, X'_k]\\
&&\qquad\qquad\qquad \hbox{\rm (since our coupling construction
is Markovian)}\\
&&\qquad\leq  \ \alpha^{k+1} B^{-N_{k-1}}
\E[h(X_{k+1}, X'_{k+1})  \mid X_k, X'_k]
\, {\bf 1}(X_k\neq X'_k)\\
&&\qquad = \ M_k \, \alpha \,
\E[h(X_{k+1}, X'_{k+1})  \mid X_k, X'_k] \, / \, h(X_k, X'_k)\cr
&&\qquad \leq  \ M_k ,
\end{eqnarray*}
by~(\ref{e9}).
Similarly, if $(X_k, X'_k)\in C\times C$, then $N_k=N_{k-1}+1$, so
assuming $X_k\neq X'_k$ (since if $X_k=X'_k$, then the result is trivial),
we have
\begin{eqnarray*}
&&\E[M_{k+1} \mid  X_0,\ldots, X_k, X'_0,\ldots, X'_k]\qquad\qquad\qquad\\
&&\qquad = \ \alpha^{k+1} B^{-N_{k-1}-1} \E[h(X_{k+1}, X'_{k+1})
\, {\bf 1}(X_{k+1}\neq X'_{k+1})  \mid  X_k. X'_k]\\
&&\qquad = \ \alpha^{k+1} B^{-N_{k-1}-1} (1-\epsilon ) (\bar{R}h)(X_k,
X'_k)\\
&&\qquad = \ M_k \alpha B^{-1} (1-\epsilon)
(\bar{R}h)(X_k, X'_k)/h(X_k, X'_k)\\
&&\qquad\leq \ M_k ,
\end{eqnarray*}
by~(\ref{e10}). Hence, $\{M_k\}$ is a supermartingale.
\end{proof}




To proceed, we note that since $B \ge 1$,
\beq
&&\P[X_k\not=X'_k, \ N_{k-1} < j]
 =  \P[X_k\not=X'_k, \ N_{k-1} \le j-1]\\
 &&\qquad \le  \P[X_k\not=X'_k, \ B^{-N_{k-1}} \ge B^{-(j-1)}]\\
 &&\qquad =  \P[\one(X_k\not=X'_k) \, B^{-N_{k-1}} \ge B^{-(j-1)}]\\
&&\qquad \le  B^{j-1} \E[\one(X_k\not=X'_k) \, B^{-N_{k-1}}]
\qquad \hbox{\rm (by Markov's inequality)}\\
&&\qquad  \le  B^{j-1} \E[\one(X_k\not=X'_k) \, B^{-N_{k-1}} \, h(X_k,X'_k)]
\qquad ({\rm since} \ h \ge 1)\\
&&\qquad  =  \alpha^{-k} B^{j-1} \E[M_k]
\qquad ({\rm by \ defn \ of} \ M_k)\\
 &&\qquad \le  \alpha^{-k} B^{j-1} \E[M_0]
\qquad ({\rm since} \ \{M_k\} \ {\rm is \ supermartingale})\\
&&\qquad  =  \alpha^{-k} B^{j-1} \E[h(X_0,X'_0)]
\qquad ({\rm by \ defn \ of} \ M_0) .
\eeq

Theorem~\ref{quantbound} now follows (in the case $n_0=1$),
by combining these two bounds with~(\ref{twoterms})
and~(\ref{couplinginequality}).

Finally, we consider the changes required if $n_0>1$.  In this case,
the main change is that we do not wish to count visits to $C \times C$
during which the joint chain could not try to couple, i.e.\ visits which
correspond to the ``filling in'' times for going back and constructing
$X_{n+1},\ldots,X_{n+n_0}$ [and similarly for $X'$] in step~2 of the
coupling construction.  Thus, we instead let $N_k$ count the number of
visits to $C \times C$, and $\{\tau_i\}$ the actual visit times, {\it
avoiding} all such ``filling in'' times.  Also, we replace $N_{k-1}$ by
$N_{k-n_0}$ in~(\ref{twoterms}) and in the definition of $M_k$.  Finally,
what is a supermartingale is not $\{M_k\}$ but rather $\{M_{t(k)}\}$,
where $t(k)$ is the latest time $\le k$ which does not correspond to a
``filling in'' time.  (Thus,  $t(k)$ will take the value $k$, unless
the joint chain visited $C \times C$ at some time between $k-n_0$
and $k-1$.) With these changes, the
proof goes through just as before.\qed

\subsection{Proof of Theorem~\ref{geomthm}}
\label{S4.5}

Here we give a direct coupling proof of Theorem~\ref{geomthm},
thereby somewhat avoiding the technicalities of e.g.\ Meyn and
Tweedie~\cite{MT} (though admittedly with a slightly weaker conclusion;
see Fact~\ref{Vunifremark}).
Our approach shall be to make use of Theorem~\ref{quantbound}.  To
begin, set $h(x,y) = \half [V(x)+V(y)]$.
Our proof will use the following technical result.

\sc{Lemma}{13}
\begin{Lemma}
\label{Vfinitelemma}\label{L14}
We may assume without loss of generality that
\be
\sup_{x \in C} V(x)  <  \infty .
\label{Vfiniteeqn}
\label{e15}
\ee
Specifically, given a small set $C$ and drift function $V$
satisfying~{\rm(\ref{minorcond})}
and~{\rm(\ref{unidriftcond})}, we can find a small set $C_0 \subseteq C$ such
that~{\rm(\ref{minorcond})}
and~{\rm(\ref{unidriftcond})}
still hold (with the same $n_0$ and $\epsilon$
and $b$, but with $\lambda$ replaced by some $\lambda_0<1$), and such
that~{\rm(\ref{Vfiniteeqn})} also holds.
\end{Lemma}


\begin{proof}
Let $\lambda$ and $b$ be as in~(\ref{unidriftcond}).
Choose $\delta$ with $0<\delta<1-\lambda$, let
$\lambda_0 = 1-\delta$, let $K = b / (1-\lambda-\delta)$, and set
$$
C_0  =  C \cap \{x\in\X: V(x) \le K\}  .
$$
Then clearly~(\ref{minorcond}) continues to hold on $C_0$, since
$C_0 \subseteq C$.
It remains to verify that~(\ref{unidriftcond}) holds
with $C$ replaced by $C_0$, and $\lambda$ replaced by $\lambda_0$.
Now,~(\ref{unidriftcond}) clearly holds
for $x \in C_0$ and $x \notin C$, by inspection.
Finally, for $x \in C \setminus C_0$, we have $V(x) \ge K$, and so
using the original drift condition~(\ref{unidriftcond}), we have
$$
(PV)(x)  \le  \lambda \, V(x) + b \, \one_C(x)
 =  (1-\delta) V(x) - (1-\lambda-\delta) V(x) + b
$$
$$
 \le  (1-\delta) V(x) - (1-\lambda-\delta) K + b
 =  (1-\delta) V(x)
 =  \lambda_0 \, V(x) ,
$$
showing that~(\ref{unidriftcond}) still holds, with $C$ replaced by $C_0$
and $\lambda$ replaced by $\lambda_0$.
\end{proof}\eject

As an aside, we note that in Lemma~\ref{Vfinitelemma}, it may
not be possible to satisfy~(\ref{Vfiniteeqn}) by instead modifying $V$
and leaving $C$ unchanged:

\sc{Proposition}{14}
\begin{Proposition}
There exists a geometrically ergodic Markov chain,
with small set $C$ and drift function $V$
satisfying~{\rm(\ref{minorcond})}
and~{\rm(\ref{unidriftcond})}, such
that there does not exist a drift function
$V_0 : \X \to [0,\infty]$ with the property that upon replacing $V$ by
$V_0$,~{\rm(\ref{minorcond})}
and~{\rm(\ref{unidriftcond})} continue to hold,
and~{\rm(\ref{Vfiniteeqn})} also holds.
\end{Proposition}

\begin{proof}
Consider the Markov chain on $\X = (0,\infty)$, defined as follows.
For $x \ge 2$, $P(x,\cdot) = \delta_{x-1}(\cdot)$, a point-mass at $x-1$.
For $1 < x < 2$, $P(x,\cdot)$ is uniform on $[1/2, \ 1]$.
For $0 < x \le 1$, $P(x,\cdot) = \half \, \lambda(\cdot) + \half \,
\delta_{h(x)}(\cdot)$, where $\lambda$
is Lebesgue measure on $(0,1)$, and $h(x) = 1 + \sqrt{ \log(1/x) }$.

For this chain, the interval $C = (0,1)$
is clearly $(1, \, 1/2, \, \lambda)$-small.  Furthermore, since
$\int_0^1 \sqrt{\log(1/x)} \, dx = \sqrt{\pi}/2 < \infty$, the return
times to $C$ have finite mean, so the chain has a stationary
distribution by standard renewal theory arguments (see
e.g.\ Asmussen~\cite{asmussen}).  In addition,
with drift function $V(x) = \max(e^x, x^{-1/2})$, we can compute $(PV)(x)$
explicitly, and verify directly
that (say) $PV(x) \le 0.8 \, V(x) + 4 \, \one_C(x)$ for all $x\in\X$,
thus verifying~(\ref{unidriftcond})
with $\lambda=0.8$ and $b=4$.  Hence, by Theorem~\ref{geomthm},
the chain is geometrically ergodic.\looseness=1

On the other hand,
suppose we had a some drift function $V_0$
satisfying~(\ref{unidriftcond}), such that
$\sup_{x \in C} V_0(x) < \infty$.  Then since $PV_0(x) = \half V_0(h(x))
+ \half \int_0^1 V_0(y) \, dy$, this would imply that
$\sup_{x \in C} V_0((h(x)) < \infty$, i.e.\ that
$V_0(h(x))$ is bounded for all $0 < x \le 1$, which would in
turn imply that $V_0$ were bounded everywhere on $\X$.
But then Fact~\ref{Vunifremark} would imply that the chain is
uniformly ergodic, which it clearly is not.  This gives a contradiction.
\end{proof}

\newcommand{\eqref}[1]{(\ref{#1})}

Thus, for the remainder of this proof, we can (and do) assume
that~\eqref{Vfiniteeqn} holds.  This, together
with~(\ref{unidriftcond}), implies that
\be
\sup_{(x,y) \in C \times C} \overline{R}h(x,y)  <  \infty ,
\label{Bfiniteeqn}\label{e16}
\ee
which in turn ensures that the quantity $B_{n_0}$
of~\eqref{Bdef} is finite.





To continue, let $d = \inf_{C^c} V$.
Then we see from Proposition~\ref{unibiprop} that the bivariate drift
condition~\eqref{bivardriftcond} will hold, {\it provided} that $d > b /
(1-\lambda) - 1$.  In that case, Theorem~\ref{geomthm} follows
immediately (in fact, in a {\it quantitative} version) by combining
Proposition~\ref{unibiprop} with Theorem~\ref{quantbound}.

However, if $d \le b / (1-\lambda) - 1$, then this argument does not
go through.  This is not merely a technicality; the condition $d > b /
(1-\lambda) - 1$ ensures that the chain is {\it aperiodic}, and without
this condition we must somehow use the assumption aperiodicity more
directly in the proof.

Our plan shall be to {\it enlarge} $C$ so that the new value of $d$
satisfies $d > b / (1-\lambda) - 1$, and to
use aperiodicity\vadjust{\eject} to
show that $C$ remains a small set (i.e., that~\eqref{minorcond}
still holds though perhaps for uncontrollably larger $n_0$ and smaller
$\epsilon>0$).  Theorem~\ref{geomthm} will then follow from
Proposition~\ref{unibiprop} and Theorem~\ref{quantbound}
as above.  (Note that we will have no direct control over the new values
of $n_0$ and $C$, which is why this approach does not provide a {\it
quantitative} convergence rate bound.)

To proceed, choose any $d' > b / (1-\lambda) - 1$, let $S = \{x\in\X; \
V(x) \le d\}$, and set $C' = C \cup S$.  This ensures that $\inf_{x\in
C'^c} V(x) \ge d' > b / (1-\lambda) - 1$.  Furthermore, since $V$
is bounded on $S$ by construction, we see that~\eqref{Vfiniteeqn} will
still hold with $C$ replaced by $C'$.  It then follows from~\eqref{Bfiniteeqn}
and~\eqref{unidriftcond} that we will still have $B_{n_0}<\infty$ even upon
replacing $C$ by $C'$.  Thus, Theorem~\ref{geomthm} will follow
from Proposition~\ref{unibiprop} and Theorem~\ref{quantbound}
if we can prove:

\sc{Lemma}{15}
\begin{Lemma}
\label{stillsmall}\label{L16}
$C'$ is a small set.
\end{Lemma}

To prove Lemma~\ref{stillsmall}, we use the notion of ``petite
set'', following~\cite{MT}.

\begin{defn}
A subset $C \subseteq \state$ is {\it petite} (or,
$(n_0,\epsilon,\nu)$-petite), relative to a small set $C$,
if there exists a positive
integer $n_0$, $\epsilon >0$, and a probability measure $\nu(\cdot)$ on
$\X$ such that
\be
\sum_{i=1}^{n_0} P^{i}(x,\cdot) \, \ge \, \epsilon \, \nu(\cdot)
\qquad x\in C  .
\label{petitecond}
\label{e17}
\ee
\end{defn}

Intuitively, the definition of petite set is like that of small
set, except that it allows the different states in $C$ to cover the
minorisation measure $\epsilon \, \nu(\cdot)$ at different times $i$.
Obviously, any small set is petite.  The converse is false in general,
as the petite set condition does not itself rule out {\it periodic}
behaviour of the chain (for example, perhaps some of the states $x\in C$
cover $\epsilon \, \nu(\cdot)$ only at odd times, and others only at even
times).  However, for an {\it aperiodic}, $\phi$-irreducible Markov chain,
we have the following result, whose proof is presented in the Appendix.

\begin{Lemma}
\label{petitelemma2}\label{L17}
(Meyn and Tweedie~\cite{MT}, Theorem~5.5.7)
For an aperiodic,\break $\phi$-irreducible Markov chain, all petite sets are
small sets.
\end{Lemma}

To make use of Lemma~\ref{petitelemma2}, we use the following.

\begin{Lemma}
\label{petitelemma1}\label{L18}
Let $C' = C \cup S$ where $S = \{x\in\X; \ V(x) \le d\}$
for some $d<\infty$, as above.  Then $C'$ is petite.
\end{Lemma}

\begin{proof}
To begin, choose $N$ large enough that $r \equiv 1 - \lambda^N d > 0$.
Let $\tau_C = \inf\{n \ge 1; \ X_n \in C\}$ be the first return time
to~$C$.  Let $Z_n = \lambda^{-n} V(X_n)$, and let $W_n =
Z_{\min(n,\tau_C)}$.  Then the drift condition~\eqref{unidriftcond}
implies that $W_n$ is a supermartingale.  Indeed, if
$\tau_C \leq n$, then
$$
\E[W_{n+1} \mid  X_0, X_1, \ldots, X_n]
 =
\E[Z_{\tau_c} \mid  X_0, X_1, \ldots, X_n]
 =
Z_{\tau_c}
 =
W_n  ,
$$\vspace*{-14pt}\eject
while if $\tau_C > n$, then $X_n \notin C$, so
using~\eqref{unidriftcond},
\begin{eqnarray*}
\E[W_{n+1} \mid  X_0, X_1, ...X_n]
 &=&
\lambda^{-(n+1)}(PV)(X_n)\\
&\leq & \lambda^{-(n+1)} \lambda V(X_n) \\
&= & \lambda^n V(X_n)\\
&= & W_n .
\end{eqnarray*}
Hence, for $x \in S$,
using Markov's inequality and the fact that $V \ge 1$,
\beq
&&\P[\tau_C \ge N \mid  X_0=x]
= \P[\lambda^{-\tau_C} \ge \lambda^{-N} \mid  X_0=x]\\
&&\qquad \le \lambda^N \, \E[\lambda^{-\tau_C} \mid  X_0=x]
\le \lambda^N \, \E[Z_{\tau_C} \mid  X_0=x]\\
&&\qquad \le \lambda^N \, \E[Z_{0} \mid  X_0=x]
= \lambda^N \, V(x)
\le \lambda^N \, d ,
\eeq
so that $\P[\tau_C < N \mid  X_0=x] \ge r$.

On the other hand, recall that $C$ is $(n_0,\epsilon,\nu(\cdot))$-small,
so that $P^{n_0}(x,\cdot) \ge \epsilon \, \nu(\cdot)$ for $x\in C$.
It follows that for $x \in S$,
$\sum_{i=1+n_0}^{N+n_0} P^i(x,\cdot) \ge r \, \epsilon \, \nu(\cdot)$.  Hence,
for $x \in S \cup C$,
$\sum_{i=n_0}^{N+n_0} P^i(x,\cdot) \ge r \, \epsilon \, \nu(\cdot)$.  This
shows that $S \cup C$ is petite.
\end{proof}e

Combining Lemmas~\ref{petitelemma1}
and \ref{petitelemma2}, we see that $C'$ must be small,
proving Lemma~\ref{stillsmall}, and hence proving
Theorem~\ref{geomthm}.\qed

\subsection{Proof of Theorem~\ref{asympconvthm}}
\label{S4.6}

Theorem~\ref{asympconvthm} does not assume the existence of any small
set $C$, so it is not clear how to make use of our coupling construction
in this case.  However, help is at hand in the form of a remarkable result
about the existence of small sets, due to Jain and Jameson \cite{jain}
(see also Orey~\cite{orey}).  We shall not prove it here; for modern
proofs see e.g.\ \cite{nummelin}, p.~16, or \cite{MT}, Theorem~5.2.2.
The key idea (see e.g.\ Meyn and Tweedie~\cite{MT}, Theorem~5.2.1) is to
extract the part of $P^{n_0}(x,\cdot)$ which is absolutely continuous with
respect to the measure $\phi$, and then to find a $C$ with $\phi(C)>0$
such that this density part is at least $\delta>0$ throughout~$C$.

\sc{thm}{18}
\begin{thm}
\label{T19}
(Jain and Jameson \cite{jain})
\label{jainthm}
Every $\phi$-irreducible Markov chain, on a state
space with countably generated $\sigma$-algebra, contains a small set
$C\subseteq\X$ with $\phi(C)>0$.  (In fact, each $B \subseteq \X$ with
$\phi(B)>0$ in turn contains a small set $C \subseteq B$ with $\phi(C)>0$.)
Furthermore, the minorisation
measure $\nu(\cdot)$ may be taken to satisfy $\nu(C)>0$.
\end{thm}

In terms of our coupling construction, if we can show that the pair
$(X_n,X'_n)$ will hit $C \times C$ infinitely often, then they will
have infinitely many opportunities to couple, with probability $\ge
\epsilon>0$ of coupling each time.  Hence, they will eventually couple
with probability~1, thus proving Theorem~\ref{asympconvthm}.

We prove this following the outline of~\cite{gentheory}.  We begin with
a lemma about return probabilities:
\sc{Lemma}{19}
\begin{Lemma}
\label{hit}\label{L20}
Consider a Markov chain on a state space $\X$, having stationary
distribution~$\pi(\cdot)$.
Suppose that for some $A \subseteq \X$, we have
$\P_x(\tau_A<\infty) > 0$ for all $x\in\X$.
Then for $\pi$-almost-every
$x\in\X$, $\P_x(\tau_A<\infty) = 1$.
\end{Lemma}

\begin{proof}
Suppose to the contrary that the conclusion does not hold, i.e.\ that
\be
\pi\big\{ x\in\X: \ \P_x(\tau_A=\infty) > 0 \big\} > 0  .
\label{contradictionassumption}
\label{e18}
\ee
Then we make the following claims (proved below):

\begin{Claim}
\label{Cl1}
Condition \eqref{contradictionassumption} implies that
there are constants $\ell, \ell_0 \in\IN$, $\delta>0$, and $B \subseteq \X$
with $\pi(B)>0$, such that
$$
\P_x \left( \tau_A=\infty,
\ \sup\{k \ge 1 ; \ X_{k \ell_0} \in B\} < \ell
\right) \ge \delta  , \qquad x\in B  .
$$
\end{Claim}

\begin{Claim}
\label{Cl2}
Let $B$, $\ell$, $\ell_0$, and $\delta$ be as in Claim~\ref{Cl1}.
Let $L = \ell \ell_0$, and let
$S = \sup\{k \ge 1; \ X_{kL} \in B\}$,
using the convention that $S=-\infty$ if the set
$\{k \ge 1; \ X_{kL} \in B\}$ is empty.
Then for all integers $1 \le r \le j$,
$$
\int_{x\in\X} \pi(dx) \, \P_x[S=r, \ X_{jL} \notin A]
 \ge  \pi(B) \, \delta .
$$
\end{Claim}

Assuming the claims, we complete the proof as follows.  We have
by stationarity that for any $j\in\IN$,
\beq
\pi(A^C) & = & \int_{x\in\X} \pi(dx) \, P^{jL}(x,A^C)
 =  \int_{x\in\X} \pi(dx) \, \P_x[X_{jL} \notin A]\\[3pt]
&\ge & \sum_{r=1}^j
\int_{x\in\X} \pi(dx) \, \P_x[S=r, \ X_{jL} \notin A]
 \ge  \sum_{r=1}^j \pi(B) \, \delta
 =  j \, \pi(B) \, \delta  .
\eeq
For $j > 1 / \pi(B) \, \delta$, this gives $\pi(A^C)>1$,
which is impossible.  This gives a contradiction, and hence
completes the proof of Lemma~\ref{hit},
subject to the proofs of Claims~\ref{Cl1} and~\ref{Cl2} below.
\end{proof}

\begin{PfC1}
By \eqref{contradictionassumption}, we can find $\delta_1$ and a subset
$B_1\subseteq\X$ with $\pi(B_1)>0$, such that
$\P_x(\tau_A<\infty) \le 1 - \delta_1$ for all $x\in B_1$.
On the other hand, since $\P_x(\tau_A<\infty) > 0$ for all $x\in\X$,
we can find $\ell_0\in\IN$ and $\delta_2>0$ and $B_2\subseteq B_1$ with
$\pi(B_2)>0$ and with $P^{\ell_0}(x,A) \ge \delta_2$ for all $x\in B_2$.

Set $\eta = \#\{k \ge 1; \ X_{k \ell_0} \in B_2\}$.  Then for any
$r\in\IN$ and $x\in\X$, we have
$\P_x(\tau_A=\infty, \ \eta = r) \le (1-\delta_2)^r$.  In particular,
$\P_x(\tau_A=\infty, \ \eta = \infty) = 0$.  Hence, for $x\in B_2$, we have
\beq
\P_x(\tau_A=\infty, \ \eta<\infty)
 &=&  1 - \P_x(\tau_A=\infty, \ \eta=\infty)
- \P_x(\tau_A<\infty)\\[3pt]
&\ge& 1 - 0 + (1-\delta_1) = \delta_1 .
\eeq
Hence, there is $\ell\in\IN$, $\delta>0$, and $B \subseteq B_2$
with $\pi(B)>0$, such that
$$
\P_x \left( \tau_A=\infty,
\ \sup\{k \ge 1 ; \ X_{k\ell_0} \in B_2\} < \ell
\right) \ge \delta  , \qquad x\in B  .
$$\eject
\noindent Finally, since $B \subseteq B_2$,
we have $\sup\{k \ge 1 ; \ X_{k\ell_0} \in B_2\}
\ge \sup\{k \ge 1 ;\break \ X_{k\ell_0} \in B\}$, thus
establishing the claim.\qed
\end{PfC1}

\begin{PfC2}
We compute using stationarity and then Claim~\ref{Cl1} that
\beq
&&\int_{x\in\X} \pi(dx) \, \P_x[S=r, \ X_{jL} \notin A]\\[6pt]
&&\qquad = \int_{x\in\X} \pi(dx) \, \int_{y\in B} P^{rL }(x,dy)
\, \P_y[S=-\infty, \ X_{(j-r)L} \notin A]\\[6pt]
&&\qquad = \int_{y\in B} \int_{x\in\X} \pi(dx) P^{rL }(x,dy)
\, \P_y[S=-\infty, \ X_{(j-r)L} \notin A]\\[6pt]
&&\qquad = \int_{y\in B} \pi(dy) \, \P_y[S=-\infty, \ X_{(j-r)L} \notin A]\\
&&\qquad \ge \int_{y\in B} \pi(dy) \, \delta
= \pi(B) \, \delta .\qed
\eeq
\end{PfC2}

To proceed, we let $C$ be a small set as in Theorem~\ref{jainthm}.
Consider again the coupling construction $\{(X_n,Y_n)\}$.  Let $G
\subseteq \X\times\X$ be the set of $(x,y)$ for which $\P_{(x,y)} \big(
\exists \, n \ge 1; \ X_n=Y_n \big) = 1$.  From the coupling construction,
we see that if $(X_0,X'_0) \equiv (x,X'_0) \in G$, then $\lim_{n\to\infty}
\P[X_n=X'_n] = 1$, so that $\lim_{n\to\infty} \|P^n(x,\cdot)
- \pi(\cdot)\| = 0$, proving Theorem~\ref{asympconvthm}.
Hence, it suffices to show that for $\pi$-a.e.\ $x\in\X$, we have
$\P[(x,X'_0) \in G] = 1$.

Let $G$ be as above, let $G_x = \{y\in\X; \ (x,y)\in G\}$ for $x\in\X$,
and let $\overline{G} = \left\{ x\in\X; \ \pi(G_x) = 1 \right\}$.  Then
Theorem~\ref{asympconvthm} follows from:

\sc{Lemma}{20}
\begin{Lemma}
$\pi(\overline{G})=1$.
\end{Lemma}

\begin{proof}
We first prove that $(\pi\times\pi)(G)=1$.
Indeed, since $\nu(C)>0$ by Theorem~\ref{jainthm},
it follows from Lemma~\ref{aperiodic} that, from any $(x,y) \in
\X\times\X$, the joint chain has positive probability of eventually
hitting $C \times C$.  It then follows by
applying Lemma~\ref{hit} to the joint chain,
that the joint chain will return to $C \times C$ with probability~1
from $(\pi\times\pi)$-a.e.\ $(x,y) \notin C \times C$.  Once
the joint chain
reaches $C \times C$, then conditional
on not coupling, the joint chain will
update from $\overline{R}$ which must be absolutely continuous with
respect to $\pi \times \pi$, and
hence (again by Lemma~\ref{hit}) will return again to $C \times C$
with probability~1.
Hence, the joint chain will repeatedly
return to $C \times C$ with
probability~1, until such time as $X_n=X'_n$.
And by the coupling construction, each time the joint chain is in $C
\times C$, it has probability
$\ge \epsilon$ of then forcing $X_n=X'_n$.
Hence, eventually we will have $X_n=X'_n$, thus proving that
$(\pi\times\pi)(G)=1$.\looseness=1\eject

Now, if we had $\pi(\overline{G})<1$, then we would have
$$
(\pi\times\pi)(G^C)  =  \int_\X \pi(dx) \pi(G_x^C)
 =  \int_{\overline{G}^C} \pi(dx) [1 - \pi(G_x)]  >  0 ,
$$
contradicting the fact that $(\pi\times\pi)(G)=1$.
\end{proof}

\section{Central Limit Theorems for Markov Chains}
\label{sec-clt}
\label{S5}

Suppose $\{X_n\}$ is a Markov chain on a state space $\X$ which is
$\phi$-irreducible and aperiodic, and has a stationary distribution
$\pi(\cdot)$.  Assume the chain begins in stationarity, i.e.\ that
$X_0 \sim \pi(\cdot)$.  Let $h:\X\to\IR$ be
some functional with finite
stationary mean $\pi(h) \equiv \int_{x\in\X} h(x) \, \pi(dx)$.

We say that $h$ satisfies a Central Limit Theorem (CLT) (or,
$\sqrt{n}$-CLT) if there is some $\sigma^2<\infty$
such that the normalised sum $n^{-1/2} \sum_{i=1}^n [h(X_i)
- \pi(h) ]$ converges weakly to a $N(0, \, \sigma^2)$ distribution.
(We allow for the special case $\sigma^2=0$, corresponding to
convergence to the constant~0.)
It then follows (see e.g.\ Chan and Geyer~\cite{changeyer}) that
\be
\sigma^2 = \lim_{n\to\infty} {1 \over n} \E\Bigg[ \Bigg(
\sum_{i=1}^n [h(X_i) - \pi(h)] \Bigg)^2 \Bigg] ,
\label{sigmaformula}
\label{e19}
\ee
and also
$\sigma^2 = \tau \, \Var_\pi(h)$, where $\tau = \sum_{k\in\IZ}
\Corr(X_0,X_k)$ is the {\it integrated autocorrelation time}.  (In the
reversible case this is also related to spectral measures; see e.g.\
\cite{kipnis}, \cite{geyerstatsci}, \cite{hybrid}.)  Clearly
$\sigma^2<\infty$ requires that $\Var_\pi(h)<\infty$, i.e.\ that
$\pi(h^2)<\infty$.

Such CLTs are helpful for understanding the {\it errors} which arise
from Monte Carlo estimation, and are thus the subject of considerable
discussion in the MCMC literature (e.g.\ \cite{geyerstatsci},
\cite{tierney}, \cite{changeyer}, \cite{MT}, \cite{RTmet},
\cite{hobertregen}, \cite{yves}, \cite{jonesclt}).

\subsection{A Negative Result}
\label{S5.1}

One might expect that CLTs always hold
when $\pi(h^2)$ is finite, but this
is false.  For example, it is shown
in~\cite{negclt} that Metropolis-Hastings algorithms whose acceptance
probabilities are too low may get so ``stuck'' that $\tau=\infty$ and
they will not have a $\sqrt{n}$-CLT.  More specifically, the following
is proved:

\sc{thm}{21}
\begin{thm}
\label{negcltthm}
Consider a reversible Markov chain, beginning in its stationary
distribution $\pi(\cdot)$, and let
$r(x) = \P[X_{n+1}=X_n \mid  X_n=x]$.  Then if
\be
\lim_{n\to\infty} \ n \ \pi\big( [h-\pi(h)]^2 \, r^n \big)  =  \infty  ,
\label{infsigma}
\label{e20}
\ee
then a $\sqrt{n}$-CLT does \un{not} hold for $h$.
\end{thm}\eject

\begin{proof}
We compute directly from~\eqref{sigmaformula} that
\beq
\sigma^2
&=& \lim_{n\to\infty} {1 \over n} \E\Bigg[ \Bigg( \sum_{i=1}^n [h(X_i)
- \pi(h) \Bigg)^2 \Bigg]\\
&\ge& \lim_{n\to\infty} {1 \over n}
\E\Bigg[ \Bigg( \sum_{i=1}^n [h(X_i)
- \pi(h) \Bigg)^2 \one(X_0=X_1=\ldots=X_n) \Bigg]\\
&=& \lim_{n\to\infty} {1 \over n} \E\Bigg[ \Bigg( n [h(X_0)
- \pi(h)] \Bigg)^2 \, r(X_0)^n \Bigg]\\
&=& \lim_{n\to\infty} n \, \pi\big( [h-\pi(h)]^2 \, r^n \big)
= \infty ,
\eeq
by~\eqref{infsigma}.  Hence, a $\sqrt{n}$-CLT cannot exist.
\end{proof}

In particular, Theorem~\ref{negcltthm} is used in~\cite{negclt}
to prove that for the independence sampler with target $\Exp(1)$
and i.i.d.\ proposals $\Exp(\lambda)$, the identity function has no
$\sqrt{n}$-CLT for any $\lambda \ge 2$.

The question then arises of what conditions on the Markov chain
transitions, and on the functional $h$, guarantee a $\sqrt{n}$-CLT
for~$h$.

\subsection{Conditions Guaranteeing CLTs}
\label{S5.2}

Here we present various positive results about the existence of CLTs.
Some, though not all, of these results are then proved in the following
two sections.

For i.i.d.\ samples, classical theory guarantees a CLT provided the
second moments are finite (e.g.\ \cite{billingsley}, Theorem~27.1;
\cite{grprobbook}, p.~110).  For {\it uniformly ergodic} chains,
an identical result exists; it is shown in Corollary~4.2(ii) of
Cogburn~\cite{cogburn} (cf.\ Theorem~5 of Tierney~\cite{tierney}) that:

\begin{thm}
\label{unifclt}
If a Markov chain with stationary distribution $\pi(\cdot)$
is uniformly ergodic, then a $\sqrt{n}$-CLT holds for
$h$ whenever $\pi(h^2)< \infty$.
\end{thm}


If a chain is just {\it geometrically ergodic} but not uniformly ergodic,
then a similar result holds under the slightly stronger assumption of a
finite $2+\delta$ moments.  That is, it is shown in Theorem~18.5.3 of
Ibragimov and Linnik~\cite{ibragimov} (see also Theorem~2 of Chan and
Geyer~\cite{changeyer}, and Theorem~2 of Hobert et al.~\cite{hobertregen})
that:

\begin{thm}
\label{geomclt}
If a Markov chain with stationary distribution $\pi(\cdot)$
is geometrically ergodic, then
a $\sqrt{n}$-CLT holds for $h$ whenever $\pi(|h|^{2+\delta})< \infty$
for some $\delta>0$.
\end{thm}

It follows, for example, that
the independence sampler example mentioned
above (which fails to have a $\sqrt{n}$-CLT, but which has finite
moments of all orders) is not geometrically ergodic.

It is shown in Corollary~3 of~\cite{hybrid} that
Theorem~\ref{geomclt} can be strengthened if the chain is {\it
reversible}:

\begin{thm}
\label{reversibleclt}
If the Markov chain is geometrically ergodic and reversible,
then a $\sqrt{n}$-CLT holds for $h$ whenever $\pi(h^2)< \infty$.
\end{thm}

%

    Comparing Theorems \ref{reversibleclt} and \ref{geomclt} leads to the following yes-or-no
    question (see \cite{yves}): if a Markov chain is geometrically ergodic,
    but not necessarily reversible, and $\pi(h^2)<\infty$, then
    does a $\sqrt{n}$-CLT necessarily exist for~$h$?  In the first
    draft of this paper, we posed that question as an Open Problem.
    However, it was recently solved by H\"aggstr\"om \cite{r36}, who
    produced a counter-example to prove the following:

    \begin{thm}
    \label{newthm26}
    There exists a (non-reversible) geometrically ergodic
    Markov chain, on a (countable) state space $\X$, and a function $h
    : \X \to \IR$, such that $\pi(h^2)<\infty$, but such that $h$ does
    not satisfy a $\sqrt{n}$-CLT (nor a CLT with any other scaling).
    \end{thm}

If $P$ is {\it reversible}, then it was proved by Kipnis and
Varadhan~\cite{kipnis} that finiteness of $\sigma^2$ is all that is
required:

\begin{thm}
\label{T26}
\label{kipnisclt}
For a $\phi$-irreducible and
aperiodic Markov chain which is reversible,
a $\sqrt{n}$-CLT holds for $h$ whenever $\sigma^2< \infty$,
where $\sigma^2$ is given
by~{\rm\eqref{sigmaformula}}.
\end{thm}

In a different direction, we have the following:

\begin{thm}
\label{VCLTthm}\label{T27}
Suppose a Markov chain is geometrically ergodic,
satisfying~\eqref{unidriftcond} for some $V:\X\to[1,\infty]$ which is
finite $\pi$-a.e.  Let
$h:\X\to\IR$ with $h^2 \le K \, V$ for some $K<\infty$.  Then
a $\sqrt{n}$-CLT holds for~$h$.
\end{thm}

Before proving some of these
results, we consider two extensions which are
straightforward mathematically, but which may be of practical importance.

\sc{Proposition}{28}
\begin{Proposition}
\label{CLTextensiona}%
The above CLT results
(i.e., Theorems~\ref{unifclt}, \ref{geomclt},
\ref{reversibleclt}, \ref{kipnisclt}, and \ref{VCLTthm})
all remain true if
instead of beginning with $X_0 \sim \pi(\cdot)$, as above, we begin
with $X_0=x$, for $\pi$-a.e.\ $x\in\X$.
\end{Proposition}

\begin{proof}
The hypotheses of the various CLT results all imply that the chain is
$\phi$-irreducible and aperiodic, with stationary distribution
$\pi(\cdot)$.  Hence, by Theorem~\ref{asympconvthm}, there is
convergence to $\pi(\cdot)$ from $\pi$-a.e.\ $x\in\X$.  For such $x$,
let $\epsilon>0$, and find $m\in\IN$ such that $\|P^m(x,\cdot) -
\pi(\cdot)\| \le \epsilon$.  It then follows from
Proposition~\ref{tvprop}(g) that we can jointly construct
copies $\{X_n\}$ and $\{X'_n\}$ of the Markov chain, with
$X_0 = x$ and $X'_0 \sim \pi(\cdot)$, such that
$$
\P[X_n = X'_n \ {\rm for \ all} \ n \ge m] \ge 1-\epsilon  .
$$\eject

\noindent But this means that for any $A \subseteq \X$,
$$
\limsup_{n\to\infty}
\Bigg| \P\Big( n^{-1/2} \sum_{i=1}^n [h(X_i) - \pi(h) ] \in A \Big)
\, - \, \P\Big( n^{-1/2}
\sum_{i=1}^n [h(X'_i) - \pi(h) ] \in A \Big) \Bigg|
 \le  \epsilon
 .
$$
Since $\epsilon>0$ is arbitrary, and since we know that
$n^{-1/2} \sum_{i=1}^n [h(X'_i) - \pi(h) ]$ converges in distribution
to $N(0,\sigma^2)$, hence
so does $n^{-1/2} \sum_{i=1}^n [h(X_i) - \pi(h) ]$.
\end{proof}

\begin{Proposition}
\label{CLTextensionb}%
The CLT Theorems~\ref{unifclt} and \ref{geomclt}
remain true if the chain is periodic of period $d \ge 2$,
provided that the $d$-step chain
$P' = P^d \big|_{\X_1}$ (as in the proof of
Corollary~\ref{periodiccor}) has all the other properties
required of $P$ in the original result
(i.e.\ $\phi$-irreducibility, and
uniform or geometric ergodicity), and
that the function $h$ still satisfies the same moment condition.
\end{Proposition}

\begin{proof}
As in the proof of Corollary~\ref{periodiccor}, let
$\Pbar$ be the $d$-step chain
defined on $\X_1 \times \ldots \times \X_d$,
and $\hbar(x_0,\ldots,x_{d-1}) = h(x_0)+\ldots+h(x_{d-1})$.
Then $\Pbar$ inherits the irreducibility and ergodicity properties
of $P'$ (formally, since $P'$ is {\it de-initialising} for $\Pbar$;
see~\cite{suffic}).  Then,
Theorem~\ref{unifclt} or \ref{geomclt} establishes a CLT for
$\Pbar$ and $\hbar$.  However, this is easily seen to be equivalent to
the corresponding CLT for
the original $P$ and $h$, thus giving the result.
\end{proof}

\begin{remark}
In particular, combining Theorem~\ref{unifclt} with
Proposition~\ref{CLTextensionb}, we see that a $\sqrt{n}$-CLT
holds for any function $h$ for any irreducible (or indecomposible)
Markov chain on a {\it finite} state space, without any assumption of
aperiodicity.  (See also the Remark following
Corollary~\ref{periodiccor} above.)
\end{remark}

\begin{remark}
We note that for periodic chains as in Proposition~\ref{CLTextensionb},
the formula~\eqref{sigmaformula} for the asymptotic variance
$\sigma^2$ continues to hold without change.
The relation $\sigma^2 = \tau \, \Var_\pi(h)$ also continues to hold, except
that now the formula for the integrated autocorrelation time
$\tau$ requires that the sum taken over ranges whose lengths are
multiples of $d$, i.e.\ the flexibly-ordered infinite sum
$\tau = \sum_{k\in\IZ} \Corr(X_0,X_k)$ must be replaced by the more
precisely limited sum
$\tau = \lim_{m,\ell\to\infty} \sum_{k=-\ell d}^{md} \Corr(X_0,X_k)$
(otherwise the sum will not converge, since now the individual terms do not
go to~0).
\end{remark}


\subsection{CLT Proofs using the Poisson Equation}
\label{S5.3}

Here we provide proofs of {\it some} of the results stated in the
previous subsection.

We begin by stating a version of the {\it martingale
central limit theorem}, which was proved independently by
Billingsley~\cite{billingsleyclt} and Ibragimov~\cite{ibragimovclt};
see e.g.~p.~375 of Durrett~\cite{durrett}.\eject

\sc{thm}{30}
\begin{thm}
\label{biclt}%
(Billingsley~\cite{billingsleyclt} and Ibragimov~\cite{ibragimovclt})
Let $\{Z_n\}$ be a stationary ergodic sequence, with
$\E[Z_n \mid  Z_1,\ldots,Z_{n-1}] = 0$ and
$\E[(Z_n)^2] < \infty$.  Then $n^{-1/2} \sum_{i=1}^n Z_i$
converges weakly to a $N(0, \, \sigma^2)$ distribution for some
$\sigma^2<\infty$.
\end{thm}

To make use of Theorem~\ref{biclt}, consider the {\it Poisson
equation}: $h - \pi(h) = g - Pg$.  A~useful result is the following
(see e.g.\ Theorem~17.4.4 of Meyn and Tweedie~\cite{MT}):

\begin{thm}
\label{poissoneqnthm}%
Let $P$ be a transition kernel for an aperiodic, $\phi$-irreducible
Markov chain on a state space $\X$,
having stationary distribution~$\pi(\cdot)$, with $X_0 \sim \pi(\cdot)$.
Let $h:\X\to\IR$ with $\pi(h^2)<\infty$, and
suppose there exists $g:\X\to\IR$ with $\pi(g^2)<\infty$ which solves
the Poisson equation, i.e.\
such that $h - \pi(h) = g - Pg$.  Then $h$ satisfies a $\sqrt{n}$-CLT.
\end{thm}

\begin{proof}
Let $Z_n = g(X_n) - Pg(X_{n-1})$.
Then $\{Z_n\}$ is stationary since $X_0 \sim \pi(\cdot)$.
Also $\{Z_n\}$ is ergodic since the Markov chain converges
asymptotically (by Theorem~\ref{asympconvthm}).
Furthermore, $\E[Z_n^2] \le 4 \, \pi(g^2) < \infty$.
Also,
\beq
\E[g(X_n)-Pg(X_{n-1}) \mid  X_0,\ldots,X_{n-1}]
 &=&  \E[g(X_n) \mid  X_{n-1}] - Pg(X_{n-1})\\
& &=  Pg(X_{n-1}) - Pg(X_{n-1})
 =  0 .
\eeq
Since $Z_1,\ldots,Z_{n-1} \in
\sigma(X_0,\ldots,X_{n-1})$, it follows
that $\E_\pi[Z_n \mid  Z_1,\ldots,Z_{n-1}] = 0$.
Hence, by Theorem~\ref{biclt},
$n^{-1/2} \sum_{i=1}^n Z_i$ converges weakly to $N(0,\sigma^2)$.  But
\beq
&&n^{-1/2} \sum_{i=1}^n[h(X_i)-\pi(h)]
 =  n^{-1/2} \sum_{i=1}^n[g(X_i)-Pg(X_i)]\\
&&\qquad  =  n^{-1/2} \sum_{i=1}^n[g(X_i)-Pg(X_{i-1})]
+n^{-1/2}Pg(X_0)-n^{-1/2}Pg(X_n)\\
&&\qquad  =  n^{-1/2}\sum_{i=1}^n Z_i
+n^{-1/2}Pg(X_0)- n^{-1/2}Pg(X_n) .
\eeq
The result follows since
$n^{-1/2} g(X_0)$ and $n^{-1/2} Pg(X_n)$ both
converge to zero in probability as $n\to\infty$.
\end{proof}

\sc{Corollary}{32}
\begin{Corollary}
\label{poissoncor}%
If $\sum_{k=0}^\infty \sqrt{\pi( (P^k[h-\pi(h)])^2)} < \infty$,
then $h$ satisfies a $\sqrt{n}$-CLT.
\end{Corollary}

\begin{proof}
Let
$$
g_k(x)  =  P^kh(x) - \pi(h)  =  P^k[h - \pi(h)](x)  ,
$$
where by convention $P^0h(x) = h(x)$, and
let $g(x) = \sum_{k=0}^\infty g_k(x)$.
Then we compute directly that
\beq
(g-Pg)(x) & =&  \sum_{k=0}^\infty g_k(x) - \sum_{k=0}^\infty Pg_k(x)
 =  \sum_{k=0}^\infty g_k(x) - \sum_{k=1}^\infty g_k(x)\\
&  = & g_0(x)
 =  P^0h(x) - \pi(h)
 =  h(x) - \pi(h) .
\eeq
Hence, the result follows from Theorem~\ref{poissoneqnthm},
{\it provided} that $\pi(g^2)<\infty$.  On the other hand,
it is known (in fact, since $\Cov(X,Y) \le \sqrt{\Var(X) \, \Var(Y)}$)
that the $L^2(\pi)$ norm satisfies the triangle inequality, so that
$$
\sqrt{\pi(g^2)}
 \le
\sum_{k=0}^\infty \sqrt{\pi(g_k^2)}
 ,
$$
so that $\pi(g^2)<\infty$ provided
$\sum_{k=0}^\infty \sqrt{\pi(g_k^2)} < \infty$.
\end{proof}

\begin{PfT25}
Let
$$
\|P\|_{L^2(\pi)} = \sup_{\pi(f)=0 \atop \pi(f^2)=1} \pi\big( (Pf)^2 \big)
= \sup_{\pi(f)=0 \atop \pi(f^2)=1} \int_{x\in\X}
\bigg( \int_{y\in\X} f(y) \, P(x,dy) \bigg)^2 \pi(dx)
$$
be the usual $L^2(\pi)$ operator norm for $P$, when
restricted to those functionals $f$ with $\pi(f)=0$
and $\pi(f^2)<\infty$.
Then it is shown in Theorem~2 of~\cite{hybrid} that reversible
chains are geometrically ergodic if and only if they satisfy
$\|P\|_{L^2(\pi)} < 1$, i.e.\ there is $\beta<1$ with
$\pi( (Pf)^2 ) \le \beta^2 \pi(f^2)$ whenever
$\pi(f)=0$ and $\pi(f^2)<\infty$.
Furthermore, reversibility implies
self-adjointness of $P$ in $L^2(\pi)$, so
that $\|P^k\|_{L^2(\pi)} = \|P\|^k_{L^2(\pi)}$,
and hence $\pi( (P^kf)^2 ) \le \beta^{2k} \pi(f^2)$.

Let $g_k = P^kh - \pi(h)$ as in the proof of
Corollary~\ref{poissoncor}.  Then this
implies that $\pi((g_k)^2) \le
\beta^{2k} \, \pi((h-\pi(h)^2)$, so that
$$
\sum_{k=0}^\infty \sqrt{ \pi(g_k^2) }
\le \sqrt{ \pi((h-\pi(h))^2) } \ \sum_{k=0}^\infty \beta^k
= \sqrt{ \pi((h-\pi(h))^2) }  /(1-\beta)
< \infty .
$$
Hence, the result follows from Corollary~\ref{poissoncor}.\qed
\end{PfT25}

\begin{PfT27}
By Fact~\ref{Vunifremark}, there is $C<\infty$ and $\rho<1$ with
$|P^nf(x) - \pi(f)| \le C V(x) \rho^n$ for $x\in\X$ and $f \le V$,
and furthermore $\pi(V)<\infty$.
Let $g_k = P^k [h-\pi(h)]$ as in the proof of
Corollary~\ref{poissoncor}.  Then
by the Cauchy-Schwartz inequality,
$(g_k)^2 = \big( P^k [h-\pi(h)] \big)^2
\le P^k\big( [h-\pi(h)]^2 \big)$.
On the other hand, since
$[h-\pi(h)]^2 \le KV$, so $[h-\pi(h)]^2/K \le V$,
we have
$(g_k)^2 \le P^k\big( [h-\pi(h)]^2 \big) \le CKV \rho^k$.
This implies that
$\pi((g_k)^2) \le CK \rho^{k} \pi(V)$, so that
\beq
&&\sum_{k=0}^\infty \sqrt{ \pi(g_k^2) }
\le \sqrt{ C K \pi((h-\pi(h))^2) }  \sum_{k=0}^\infty \rho^{k/2}\\
&&\qquad = \sqrt{ C K \pi((h-\pi(h))^2) }  / (1-\sqrt{\rho})
< \infty
 .
\eeq
Hence, the result again
follows from Corollary~\ref{poissoncor}.\qed
\end{PfT27}
\eject

\subsection{Proof of Theorem~\ref{geomclt} using Regenerations}
\label{regensec}
\label{S5.4}

Here we use {\it regeneration theory} to give a reasonably direct
proof of Theorem~\ref{geomclt}, following the outline of Hobert
et al.~\cite{hobertregen}, thereby avoiding the technicalities of the
original proof of Ibragimov and Linnik~\cite{ibragimov}.

We begin by noting from Fact~\ref{Vunifremark} that since the
chain is geometrically ergodic, there is a small set $C$ and a drift
function $V$ satisfying~\eqref{minorcond} and~\eqref{unidriftcond}.

In terms of this, we consider a {\it regeneration construction}
for the chain (cf.\ \cite{athreya}, \cite{asmussen}, \cite{mykland},
\cite{hobertregen}).  This is very similar to the coupling construction
presented in Section~\ref{S4}, except now just for a {\it single}
chain $\{X_n\}$.  Thus, in the coupling construction we omit option~1,
and merely update the single chain.  More formally, given $X_n$,
we proceed as follows.  If $X_n \notin C$, then we simply choose
$X_{n+1} \sim P(X_n,\cdot)$.  Otherwise, if $X_n \in C$, then with
probability $\epsilon$ we choose $X_{n+n_0} \sim \nu(\cdot)$, while
with probability $1-\epsilon$ we choose $X_{n+n_0} \sim R(X_n,\cdot)$.
[If $n_0>1$, we then fill in the missing values
$X_{n+1},\ldots,X_{n+n_0-1}$ as usual.]

We let $T_1,T_2,\ldots$ be the {\it regeneration times}, i.e.\ the times
such that $X_{T_i} \sim \nu(\cdot)$ as above.  Thus, the regeneration
times occur with probability $\epsilon$ precisely $n_0$ iterations after
each time the chain enters $C$ (not counting those entries of $C$ which
are within $n_0$ of a previous regeneration attempt).

The benefit of regeneration times is that they break up sums like
$\sum_{i=0}^n [h(X_i)-\pi(h)]$ into sums over {\it tours}, each of the
form $\sum_{i=T_j}^{T_{j+1}-1} [h(X_i)-\pi(h)]$.  Furthermore, since each
subsequent tour begins from the same fixed distribution $\nu(\cdot)$, we
see that {\sl the different tours, after the first one, are independent
and identically distributed (i.i.d.)}.

More specifically, let $T_0=0$, and let $r(n) = \sup\{i \ge 0; \ T_i \le n\}$.
Then
\be
\sum_{i=1}^n [h(X_i)-\pi(h)]
 =  \sum_{j=1}^{r(n)} \sum_{i=T_j}^{T_{j+1}-1} [h(X_i)-\pi(h)]
+ E(n) ,
\label{toursum}\label{e21}
\ee
where $E(n)$ is an error term which collects the terms corresponding
to the incomplete final tour $X_{T_{r(n)+1}},\ldots,X_n$,
and also the first tour
$X_0,\ldots,X_{T_1-1}$.

Now, the tours $\{\{X_{T_j}, X_{T_j+1}, \ldots ,X_{T_{j+1}-1}\}, \
j=1, 2, \ldots \}$ are independent and identically distributed.
Moreover, elementary renewal theory (see for example \cite{asmussen})
ensures that $r(n)/n \to \epsilon \pi (C)$ in probability.  Hence, the
classical central limit theorem (see e.g.\ \cite{billingsley},
Theorem~27.1; or \cite{grprobbook}, p.~110) will prove
Theorem~\ref{geomclt}, {\it provided} that each term has finite second
moment, and that the error term $E(n)$ can be neglected.

To continue, we note that geometric ergodicity
implies (as in the proof
of Lemma~\ref{petitelemma1}) exponential tails on the return times
to~$C$.  It then follows (cf.\ Theorem~2.5 of \cite{tuominen}) that
there is $\beta>1$ with
\be
\E_\pi[\beta^{T_1}]<\infty  ,
\quad {\rm and} \quad
\E[ \beta^{T_{j+1} - T_j} ] < \infty .
\label{finitereturnmomentseqn}
\label{e22}
\ee
(This also follows from Theorem 15.0.1 of \cite{MT}, together with a
simple argument using probability generating functions.)

Now, it seems intuitively
clear that $E(n)$ is $O_p(1)$ as $n\to\infty$,
so when multiplied by $n^{-1/2}$, it will not contribute to the limit.
Formally, this follows from~\eqref{finitereturnmomentseqn}, which
implies by standard renewal theory that $E(n)$ has a limiting
distribution as $n \to \infty$, which in turn implies that $E(n)$ is
$O_p(1)$ as $n\to\infty$.  Thus, the term $E(n)$ can be neglected
without affecting the result.

Hence, it remains only to
prove the finite second moments of each term
in~\eqref{toursum}.  Recalling
that each tour begins in the distribution
$\nu(\cdot)$, we see that the proof of Theorem~\ref{geomclt} is
completed by the following lemma:

\sc{Lemma}{33}
\begin{Lemma}
$\int_{x\in\X} \nu(dx) \,
\E \big[ \big( \sum_{i=0}^{T_{1}-1} [h(X_i)-\pi(h)] \big)^2
\, \bigm| \, X_{0}=x \big]
\, < \, \infty$.
\end{Lemma}

\begin{proof1}
Note that
$$
\pi(\cdot)  =  \int_{x\in\X} \pi(dx) \, P(x,\cdot)
 \ge
\int_{x\in C} \pi(dx) \, P(x,\cdot)
 \ge  \pi(C) \, \epsilon \, \nu(\cdot)  ,
$$
so that $\nu(dx) \le \pi(dx) /\pi(C) \, \epsilon$.
Hence, it suffices to prove the lemma
with $\nu(dx)$ replaced by $\pi(dx)$
i.e.\ under the assumption that $X_0 \sim \pi(\cdot)$.

For notational simplicity, set $H_i = h(X_i)-\pi(h)$, and
$\E_\pi[\cdots] = \int_{x\in\X} \E[\cdots \mid  X_{0}=x] \, \pi(dx)$.
Note that
$( \sum_{i=0}^{T_{1}-1} [h(X_i)-\pi(h)])^2
= ( \sum_{i=0}^\infty \one_{i < T_{1}} H_i )^2$.
Hence, by Cauchy-Schwartz,
\be
\label{e23}
\E_\pi\Bigg[ \Bigg( \sum_{i=0}^{T_{1}-1}
[h(X_i)-\pi(h)] \Bigg)^2 \Bigg]
 \le
\Bigg( \sum_{i=0}^\infty
\sqrt{ \E_\pi\big[ \one_{i < T_{1}} H_i^2 \big] } \ \Bigg)^2 .
\label{CSeqn}
\ee

To continue, let $p=1+2/\delta$ and $q=1+\delta/2$, so
that $1/p+1/q=1$.
Then by {\it H\"older's inequality} (e.g.~\cite{billingsley}, p.~80),
\be
\E_\pi[ \one_{i < T_{1}} H_i^2 ]
 \le
\E_\pi[ \one_{i < T_{1}} ]^{1/p}
\, \E_\pi[ |H_i|^{2q} ]^{1/q} .
\label{holdereqn}
\label{e24}
\ee

Now, since $X_0\sim\pi(\cdot)$, therefore
$\E_\pi[ |H_i|^{2q} ] \equiv K$ is a constant, independent of
$i$, which is finite since $\pi(|h|^{2+\delta})<\infty$.

Also, using~\eqref{finitereturnmomentseqn},
Markov's inequality then gives that
$\E_\pi[ \one_{0 \le i < T_{1}} ] \le\break \E_\pi[ \one_{\beta^{T_1} >
\beta^i} ] \le \beta^{-i} \E_\pi[\beta^{T_1}]$.
Hence, combining~\eqref{CSeqn} and~\eqref{holdereqn}, we obtain that
\beq
&&\E_\pi\Bigg[ \Bigg( \sum_{i=0}^{T_{1}-1}
[h(X_i)-\pi(h)] \Bigg)^2 \Bigg]
 \le
\Bigg( \sum_{i=0}^\infty \sqrt{
\E_\pi[ \one_{i < T_{1}} ]^{1/p} \, \E_\pi
[ |H_i|^{2q} ]^{1/q}}\Bigg)^2\\[6pt]
&&\qquad \le \Bigg(
K^{1/2q} \sum_{i=0}^\infty \sqrt{ (\beta^{-i}
\E_\pi[\beta^{T_1}])^{1/p} }
\Bigg)^2
= \Bigg(
K^{1/2q} \E_\pi[\beta^{T_1}]^{1/2p} \sum_{i=0}^\infty \beta^{-i/2}
\Bigg)^2\\[6pt]
&&\qquad = \big(
K^{1/2q} \E_\pi[\beta^{T_1}]^{1/2p} \, / \, (1-\beta^{-1/2})
\big)^2
< \infty .\qed
\eeq
\end{proof1}

It appears at first glance that Theorem~\ref{unifclt} could be
proved by similar regeneration
arguments.  However, we have been unable
to do so.\eject

\begin{opennew}
\label{O3}
Can Theorem~\ref{unifclt} be proved by direct regeneration arguments,
similar to the above proof of Theorem~\ref{geomclt}?
\end{opennew}

\section{Optimal Scaling and Weak Convergence}
\label{sec-optimal}
\label{S6}

Finally, we briefly discuss another
application of probability theory to
MCMC, namely the {\it optimal scaling} problem.  Our presentation here
is quite brief; for further details see the review
article~\cite{statsci}.

Let $\pi_u: {\bf R}^d \to [0,\infty)$ be a continuous $d$-dimensional
density ($d$ large).  Consider running a Metropolis-Hastings algorithm
for $\pi_u$.  The optimal scaling problem concerns the question of how
we should choose the proposal distribution for this algorithm.

For concreteness, consider either the random-walk Metropolis (RWM)
algorithm with proposal distribution given by $Q(x,\cdot) = N(x,\sigma^2
I_d)$, or the Langevin algorithm with proposal distribution given
by $Q(x,\cdot) = N(x + {\sigma^2 \over 2} \nabla \log \pi_u(x) , \
\sigma^2 I_d)$.  In either case, the question becomes, how should we
choose $\sigma^2$?

If $\sigma^2$ is chosen to be too small, then by continuity the resulting
Markov chain will nearly always accept its proposed value.  However, the
proposed value will usually be extremely close to the chain's previous
state, so that the chain will move extremely slowly, leading to a very
high acceptance rate, but very poor performance.  On the other hand,
if $\sigma^2$ is chosen to be too large, then the proposed values will
usually be very far from the current state.  Unless the chain gets very
``lucky'', then those proposed values will usually be rejected, so that
the chain will tend to get ``stuck'' at the same state for large periods
of time.  This will lead to a very low acceptance rate, and again a very
poorly performing algorithm.  We conclude that proposal scalings satisfy a
{\it Goldilocks Principle}:  The choice of the proposal scaling $\sigma^2$
should be ``just right'', neither too small nor too large.

To prove theorems about this, assume for now that
\be
\pi_u(\bx)  =  \prod _{i=1}^d f(x_i)  ,
\label{iidcomponents}
\label{e25}
\ee
i.e.\ that the density $\pi_u$ factors
into i.i.d.\ components, each with
(smooth) density $f$.  (This assumption is obviously very restrictive, and
is uninteresting in practice since then
each coordinate can be simulated
separately.  However, it does
allow us to develop some interesting theory, which may approximately
apply in other cases as well.)  Also, assume that chain begins
in stationarity, i.e.\ that $X_0 \sim \pi(\cdot)$.

\subsection{The Random Walk Metropolis (RWM) Case}
\label{S6.1}

For RWM, let $I=\E[((\log f(Z))')^2]$
where $Z \sim f(z) \, dz$.  Then it
turns out, essentially, that under the
assumption~\eqref{iidcomponents},
as $d\to\infty$ it is optimal to choose $\sigma^2 \doteq (2.38)^2 /
I d$, leading to an asymptotic acceptance rate $\doteq 0.234$.\eject

More precisely, set the
proposal variance to be $\sigma_d^2=\ell^2/d$, where
$\ell>0$ is to be chosen later.  Let
$\{X_n\}$ be the Random Walk Metropolis
algorithm for $\pi(\cdot)$ on $\IR^d$ with proposal variance
$\sigma_d^2$.
Also, let $\{N(t)\}_{t \ge 0}$ be a Poisson
process with rate $d$ which is
independent of $\{X_n\}$.  Finally, let
$$
Z_t^d = X_{N(t)}^{(1)}  ,
\qquad t \ge 0  .
$$
Thus, $\{Z_t^d\}_{t \ge 0}$ follows the first component of $\{X_n\}$,
with time speeded up by a factor of~$d$.

Then it is proved in~\cite{RGG}
(see also~\cite{statsci}), using the theory
from Ethier and Kurtz~\cite{ethier}, that as $d\to\infty$,
the process $\{Z_t^d\}_{t \ge 0}$ converges weakly to a diffusion
process $\{Z_t\}_{t \ge 0}$ which satisfies the following stochastic
differential equation:
$$
dZ_t  =  h(\ell )^{1/2} \, dB_t
+ \half \, h(\ell) \, \nabla \log \pi_u (Z_t) \ dt
 .
$$
Here
$$
h(\ell )  =  2 \, \ell^2 \, \Phi\left(-{\sqrt{I}\ell \over 2}\right)
$$
corresponds to the {\it speed} of the limiting diffusion, where
$\Phi(x) = {1 \over \sqrt{2\pi}} \int_{-\infty}^x e^{-s^2/2} ds$
is the cdf of a standard normal distribution.

We then compute numerically that the choice $\ell = \hat\ell \doteq 2.38 /
\sqrt{I}$ maximises the above speed function $h(\ell)$, and thus must be
the choice leading to optimally fast mixing (at least, as $d\to\infty$).
Furthermore, it is also proved in~\cite{RGG} that the asymptotic (i.e.,
expected value with respect to the stationary distribution) acceptance
rate of the algorithm is given by the formula $A(\ell) = 2 \, \Phi
\left(-{\sqrt{I}\ell \over 2}\right)$, and we compute that $A(\hat\ell)
\doteq 0.234$, thus giving the optimal asymptotic acceptance rate.

\subsection{The Langevin Algorithm Case}
\label{S6.2}

In the Langevin case, let $J=\E[( 5((\log f(Z))'''))^2 - 3((\log
f(Z))'')^3)/48]$ where again $Z \sim f(z) \, dz$.  Then it turns out,
essentially, that assuming~\eqref{iidcomponents}, it is optimal as
$d\to\infty$ to choose $\sigma^2 \doteq (0.825)^2 / J^{1/2} d^{1/3}$,
leading to an asymptotic acceptance rate $\doteq 0.574$.

More precisely, set $\sigma_d^2=\ell^2/d^{1/3}$, let $\{X_n\}$ be the
Langevin Algorithm for $\pi(\cdot)$ on $\IR^d$
with proposal variance $\sigma_d^2$, let $\{N(t)\}_{t \ge 0}$
be a Poisson process with rate $d^{1/3}$ which is independent of
$\{X_n\}$, and let
$$
Z_t^d = X_{N(t)}^{(1)}  ,
$$
so that $\{Z_t^d\}_{t \ge 0}$
follows the first component of $\{X_n\}$,
with time speeded up by a factor of $d^{1/3}$.
Then it is proved in~\cite{lang} (see also~\cite{statsci}) that
as $d\to\infty$,
the process $\{Z_t^d\}_{t \ge 0}$ converges weakly to a diffusion
process $\{Z_t\}_{t \ge 0}$ which satisfies the following stochastic
differential equation:
$$
dZ_t  =  g(\ell )^{1/2} \, dB_t
+ \half \, g(\ell) \, \nabla \log \pi_u (Z_t) \ dt
 .
$$
Here
$$
g(\ell )  =  2 \, \ell^2 \, \Phi\left(-J\ell^3\right)
$$
represents the speed of the limiting diffusion.
We then compute numerically that
the choice $\ell = \hat\ell = 0.825/\sqrt{J}$ maximises $g(\ell)$,
and thus must be
the choice leading to optimally fast mixing (at least, as $d\to\infty$).
Furthermore, it is proved in~\cite{lang} that
the asymptotic acceptance rate satisfies
$A(\hat\ell) = 2 \, \Phi(-J\hat\ell^3) \doteq 0.574$,
thus giving the optimal asymptotic acceptance rate for the Langevin case.

\subsection{Discussion of Optimal Scaling}
\label{S6.3}

The above results show that for either the RWM or the Langevin algorithm,
under the assumption~\eqref{iidcomponents}, we can determine the optimal
proposal scaling just in terms of universally optimal asymptotic
acceptance rates (0.234 for RWM, 0.574 for Langevin).  Such results
are straightforward to apply in practice,
since it is trivial for a computer to monitor the acceptance rate of
the algorithm, and the user can modify $\sigma^2$ appropriately to
achieve appropriate acceptance rates.
Thus, these optimal scaling rates are often used in applied contexts
(see e.g.\ M\o ller et al.~\cite{molsyvwaa}).
(It may even be possible for
the computer to adaptively modify $\sigma^2$ to achieve the appropriate
acceptance rates; see~\cite{yvesadapt} and references therein. However
it is important to recognise that adaptive strategies can violate the
stationarity of $\pi $ so they have to be carefully implemented; see
for example \cite{robsah}.)

The above results also describe the {\it computational complexity}
of these algorithms.  Specifically, they say that as $d\to\infty$, the
efficiency of RWM algorithms scales like $d^{-1}$, so its computational
complexity is $O(d)$.  Similarly, the efficiency of Langevin algorithms
scales like $d^{-1/3}$, so its computational complexity is $O(d^{1/3})$
which is much lower order (i.e.\ better).

We note that for reasonable efficiency, we do not need the acceptance
rate to be {\it exactly} 0.234 (or 0.574), just fairly close.  Also,
the dimension doesn't have to be {\it too} large before asymptotics
approximately kick in; often 0.234 is approximately optimal in dimensions
as low as~5 or~10.  For further discussion of these issues, see the
review article~\cite{statsci}.

Now, the above results are only proved under the strong
assumption~\eqref{iidcomponents}.  It is natural to ask what
happens if this assumption is not satisfied.  In that case, there
are various extensions of the optimal-scaling results to cases
of inhomogeneously-scaled components of the form $\pi_u(\bx) =
\prod_{i=1}^d C_i \, f(C_i x_i)$\break (see~\cite{statsci}), to the discrete
hypercube~\cite{hypercube}, and to finite-range homogeneous Markov random
fields~\cite{MRF}; in particular, the optimal acceptance rate remains
0.234 (under appropriate assumptions) in all of these cases.  On the other
hand, surprising behaviour can result if we do not start in stationarity,
i.e.\ if the assumption $X_0 \sim \pi(\cdot)$ is violated and the chain
instead begins way out in the tails of $\pi(\cdot)$; see~\cite{determ}.
The true level of generality of these optimal scaling results is currently
unknown, though investigations are ongoing~\cite{mylene}.  In general
this is an open problem:\eject

\begin{opennew}
Determine the extent to which the
above optimal scaling results continue
to apply, even when assumption~\eqref{iidcomponents} is violated.
\end{opennew}

\section*{APPENDIX: Proof of
Lemma~\ref{petitelemma2}}

Lemma~\ref{petitelemma2} above states (Meyn and Tweedie~\cite{MT},
Theorem~5.5.7) that for an aperiodic, $\phi$-irreducible Markov chain,
all petite sets are small sets.

To prove this, we require a lemma related to aperiodicity:

\sc{Lemma}{34}
\begin{Lemma}
\label{aperiodic}
Consider an aperiodic Markov chain on a state space
$\X$, with stationary distribution $\pi(\cdot)$.
Let $\nu(\cdot)$ be any probability measure on $\X$.
Assume that $\nu(\cdot) \ll \pi(\cdot)$, and that for all $x\in\X$,
there is $n=n(x)\in\IN$ and $\delta=\delta(x)>0$ such that
$P^n(x,\cdot) \ge \delta \nu(\cdot)$
(for example, this always holds if
$\nu(\cdot)$ is a minorisation
measure for a small or petite set which
is reachable from all states).
Let $T = \{ n \ge 1; \ \exists \delta_n > 0 \ s.t.\
\int \nu(dx) \, P^n(x,\cdot) \ge \delta_n \nu(\cdot) \}$, and assume
that $T$ is non-empty.  Then there is
$n_*\in\IN$ with $T \supseteq \{n_*,n_*+1,n_*+2,\ldots\}$.
\end{Lemma}

\begin{proof}
Since $P^{(n(x))}(x,\cdot) \ge
\delta(x) \, \nu(\cdot)$ for all $x\in\X$,
it follows that $T$ is non-empty.

Now, if $n,m \in T$, then since $\int_{x\in\X} \nu(dx) \,
P^{n+m}(x,\cdot) = \int_{x \in \X} \int_{y\in\X} \nu(dx)\times\break
P^n(x,dy) P^m(y,\cdot) \ge \int_{y\in\X} \delta_n \nu(dy) P^m(y,\cdot)
\ge \delta_n \delta_m \nu(\cdot)$, we see that $T$ is {\it additive},
i.e.\ if $n,m\in T$ then $n+m\in T$.

We shall prove below that $gcd(T)=1$.
It is then a standard and easy
fact (e.g.\ \cite{billingsley}, p.~541; or \cite{grprobbook}, p.~77)
that if $T$ is non-empty and additive, and $gcd(T)=1$, then
there is $n_*\in\IN$ such that $T \supseteq \{n_*, n_*+1,n_*+2,
\ldots\}$, as claimed.

We now proceed to prove that $gcd(T)=1$.  Indeed, suppose to the contrary
that $gcd(T)=d>1$.  We will derive a contradiction.

For $1 \le i \le d$, let
$$
\X_i  =  \{x\in\X; \ \exists \ell\in\IN \ {\rm and} \ \delta>0
\ {\rm s.t.}\ P^{\ell d - i}(x,\cdot) \ge \delta \nu(\cdot) \}
 .
$$
Then $\bigcup_{i=1}^d \X_i = \X$ by assumption.  Now, let
$$
S  =  \bigcup_{i \not= j} ( \X_i \cap \X_j )
 ,
$$
let
$$
\overline{S}  =  S \cup \{x\in\X; \ \exists m\in\IN \ s.t.\ P^m(x,S)>0\}
 ,
$$
and let
$$
\X'_i  =  \X_i \setminus \overline{S}  .
$$

Then $\X_1,\X_2,\ldots, \X_d$
are disjoint by construction (since we have removed $S$).
Also if $x\in\X'_i$, then $P(x,\overline{S})=0$,
so that $P(x,\bigcup_{j=1}^d \X'_j)=1$ by construction.
In fact we must have $P(x,\X'_{i+1})=1$ in the case $i<d$ (with
$P(x,\X'_{1})=1$ for $i=d$), for if not then $x$ would be in two
different $\X'_j$ at once, contradicting their disjointedness.\eject

We claim that for all $m \ge 0$, $\nu P^m (\X_i \cap \X_j)=0$ whenever
$i \not= j$.  Indeed, if we had $\nu P^m (\X_i \cap \X_j)>0$ for some $i
\not= j$, then there would be $S' \subseteq \X$, $\ell_1,\ell_2\in\IN$,
and $\delta>0$ such that for all $x\in S'$, $P^{\ell_1 d + i}(x,\cdot)
\ge \nu(\cdot)$ and $P^{\ell_2 > d + j}(x,\cdot) \ge \nu(\cdot)$, implying
that $\ell_1 d + i + m \in T$ and $\ell_2 d + j + m \in T$, contradicting
the fact that $\gcd(T)=d$.

It then follows (by sub-additivity of measures) that
$\nu(\overline{S})=0$.  Therefore, $\nu(\bigcup_{i=1}^d \X'_i) =
\nu(\bigcup_{i=1}^d \X_i) = \nu(\X) = 1$.  Since $\nu \ll \pi$, we must
have $\pi(\bigcup_{i=1}^d \X'_i) > 0$.

We conclude from all of this that $\X'_1,\ldots,\X'_d$ are subsets of
positive $\pi$-measure, with respect to which the Markov chain is periodic
(of period $d$), contradicting the assumption of aperiodicity.
\end{proof}

\begin{PfL17}
Let $R$ be $(n_0,\epsilon,\nu(\cdot))$-petite, so that
$\sum_{i=1}^{n_0} P^i(x,\cdot) \ge \epsilon \nu(\cdot)$ for all $x \in R$.
Let $T$ be as in Lemma~\ref{aperiodic}.
Then $\sum_{i=1}^{n_0} \int_{x\in\X} \nu(dx) \, P^i(x,\cdot)
\ge \epsilon \nu(\cdot)$, so we must have $i\in T$ for some
$1 \le i \le n_0$, so that $T$ is non-empty.
Hence, from Lemma~\ref{aperiodic}, we can find $n_*$ and
$\delta_n>0$ such that
$\int \nu(dx) P^n(x,\cdot) > \delta_n \nu(\cdot)$ for all $n \ge n_*$.
Let $r  =  \min\Big\{ \delta_n ; \ n_* \le n \le n_*+n_0-1 \Big\}$,
and set $N = n_* + n_0$.  Then for $x \in R$,
\beq
P^N(x,\cdot) &\ge&
\sum_{i=1}^{n_0} \int_{y\in\X} P^{N-i}(x,dy) P^{i}(y,\cdot)\\
&\ge& \sum_{i=1}^{n_0} \int_{y\in R} r \nu(dy) P^{i}(y,\cdot)\\
&\&ge \int_{y\in R} r \nu(dy) \epsilon \nu(\cdot)
= r \epsilon \nu(\cdot) .
\eeq
Thus, $R$ is $(N,r\epsilon,\nu(\cdot))$-small.\qed
\end{PfL17}

%
%

\section*{Acknowledgements}
The authors are sincerely grateful to Kun Zhang and to an
anonymous referee for their extremely careful readings of this
manuscript and their many insightful comments
which lead to numerous improvements.\pagebreak

\end{document}